\documentclass[10pt,a4paper,titlepage]{article}
\usepackage[mathscr,mathcal]{eucal}
\usepackage{enumerate}
\usepackage{colortbl}
\usepackage{amsfonts}
\usepackage{amssymb}
\usepackage{amsmath}
\newtheorem{theorem}{Theorem}[section]
\newtheorem{definition}[theorem]{Definition}
\newtheorem{ypoth}[theorem]{Assumptions}

\newtheorem{lemma}[theorem]{Lemma}

\newtheorem{corollary}[theorem]{Corollary}
\newtheorem{proposition}[theorem]{Proposition}
\newtheorem{thedef}[theorem]{Theorem-Definition}

\newcommand{\el}{\lambda}
\newcommand{\eb}{\beta}
\newtheorem{px}[theorem]{Example}

\newcommand{\ie}{\emph{i.e., }}

\begin{document}

\begin{center}
\textbf{Rouquier blocks of the cyclotomic Hecke algebras of $G(de,e,r)$}\\
$ $\\
Maria Chlouveraki\\
{\small Ecole Polytechnique F\'ed\'erale de Lausanne}\\
{\small \emph{e-mail : maria.chlouveraki@epfl.ch}}
\end{center}
ABSTRACT. The ``Rouquier blocks'' of the cyclotomic Hecke algebras, introduced by Rouquier, 
are a substitute for the ``families of characters'', defined by Lusztig for Weyl groups,  which can be applied to all complex reflection groups. In this article, we determine them for the cyclotomic Hecke algebras of the groups of the infinite series,
$G(de,e,r)$, thus completing their calculation for all complex reflection groups.


\section *{Introduction}

Until recently, the lack of Kazhdan-Lusztig bases for the non-Coxeter complex reflection groups did not allow the generalization of the notion of ``families of characters'' from Weyl groups to all complex reflection groups. However, thanks to the results of  Gyoja $\cite{Gy}$ and Rouquier
$\cite{Rou}$, we have obtained a substitute for the
families of characters which can be applied to all complex
reflection groups. In particular, Rouquier has proved that the families
of characters of a Weyl group $W$ coincide with the ``Rouquier blocks'' of the Iwahori-Hecke algebra of $W$, \ie its blocks over a suitable coefficient ring. This definition generalizes to all complex reflection groups and we are grateful for this for the following reasons:

On one hand, since the families of characters of a Weyl group play an essential role in the
definition of the families of unipotent characters of the
corresponding finite reductive group (cf.~\cite{Lu1}), 
the families of characters of the cyclotomic Hecke algebras could play a key
role in the organization of families of unipotent characters in general.
On the other hand, for some (non-Coxeter) complex reflection groups 
$W$, we have data which seem to indicate that behind
the group $W$, there exists another mysterious object - the
\emph{Spets} (cf. \cite{BMM2}, \cite{Ma3}) - that could play the role
of the ``series of finite reductive groups of Weyl group $W$''.

In \cite{BK}, Brou\'{e} and Kim presented an algorithm for the determination of the Rouquier blocks of the cyclotomic Hecke algebras of the groups $G(d,1,r)$. Using the generalization of some classic results, known as ``Clifford theory'',  they were able to obtain the Rouquier blocks for $G(d,d,r)$ from those of $G(d,1,r)$. Later, Kim \cite{Kim} generalized the methods used in \cite{BK} in order to obtain the Rouquier blocks of the cyclotomic Hecke algebras of $G(de,e,r)$ from those of $G(de,1,r)$.

As far as the exceptional complex reflection groups are concerned,  some special cases were treated  by Malle and Rouquier in \cite{MaRo}. Finally, in \cite{Chlou}, we gave the complete classification of the  Rouquier blocks of the cyclotomic Hecke algebras for all exceptional complex reflection groups. 

However, recently it was realized that the algorithm of \cite{BK} for $G(d,1,r)$ does not work, unless $d$ is a power of a prime number. 
In \cite{Chlou2}, we give the correct algorithm, which is more complicated than the one of \cite{BK}. 
Now, it remains to recalculate the Rouquier blocks of the cyclotomic Hecke algebras of $G(de,e,r)$, in order to complete the determination of the Rouquier blocks for all complex reflection groups.  

Using the same idea as in \cite{Kim}, we apply ``Clifford theory'' in order to obtain the Rouquier blocks for $G(de,e,r)$ from those of $G(de,1,r)$. However, we point out that there is one case where this is not possible, that is, when $r=2$ and $e$ is even. In that case, we apply the same methods as in \cite{Chlou} in order to determine the Rouquier blocks of the cyclotomic Hecke algebras of $G(de,2,2)$, and then ``Clifford theory'' in order to obtain the Rouquier blocks  for $G(de,e,2)$.

Finally, to every irreducible character of a cyclotomic Hecke algebra of a complex reflection group we can attach integers $a$ and $A$, like Lusztig has done for Weyl groups. In \cite{Lu2}, Lusztig shows that these integers are constant on families. Here, we complete the proof that $a$ and $A$ are constant on the Rouquier blocks of the cyclotomic Hecke algebras of all irreducible complex reflection groups, having already shown it for the exceptional ones (cf. \cite{ChDeg}) and $G(d,1,r)$ (cf. \cite{Chlou2}).

\section {Blocks of symmetric algebras}

All the results of this section are presented here for the convenience of the reader.
Their proofs can be found in the second chapter of \cite{Chlou}.

\subsection {Generalities on blocks}

Let us assume that
$\mathcal{O}$ is a commutative integral domain with field of fractions
  $F$ and  $A$ is an $\mathcal{O}$-algebra, free and finitely generated as an
  $\mathcal{O}$-module.

\begin{definition}\label{blocks}
The block-idempotents (blocks) of $A$ are the primitive idempotents
of $ZA$.
\end{definition}

Let $K$ be a field extension of $F$. Suppose that the $K$-algebra $KA:=K \otimes_\mathcal{O}A$ is
semisimple. Then there exists a bijection between the set $\mathrm{Irr}(KA)$  of irreducible characters of  $KA$ and the set $\mathrm{Bl}(KA)$ of blocks of $KA$:
$$\begin{array}{ccc}
   \mathrm{Irr}(KA) & \leftrightarrow & \mathrm{Bl}(KA) \\
    \chi &\mapsto &e_\chi.
  \end{array}$$

The following theorem establishes a relation between the blocks of the algebra $A$ and the blocks of $KA$.

\begin{theorem}\label{minimality of blocks}\
There exists a unique partition $\mathrm{Bl}(A)$ of
  $\mathrm{Irr}(KA)$ such that
  \begin{enumerate}[(1)]
    \item For all $B \in \mathrm{Bl}(A)$, the idempotent
    $e_B:=\sum_{\chi \in B}e_\chi$ is a block of $A$.
    \item For
    every central idempotent $e$ of $A$, there exists a subset
    $\mathrm{Bl}(A,e)$ of $\mathrm{Bl}(A)$ such that
    $$e=\sum_{B \in \mathrm{Bl}(A,e)}e_B.$$
      \end{enumerate}
    In particular, the set $\{e_B\}_{B \in \mathrm{Bl}(A)}$ is the set of all the blocks of $A$.
\end{theorem}
If $\chi \in B$ for some $B \in \mathrm{Bl}(A)$, we say that
  ``$\chi$ belongs to the block $e_B$''.

\subsection {Symmetric algebras}

From now on, we make the following assumptions

\begin{ypoth}\label{properties of the ring}\
\begin{description}
  \item[(int)] The ring $\mathcal{O}$ is a Noetherian and integrally
  closed domain with field of fractions $F$ and $A$ is an
  $\mathcal{O}$-algebra which is free and finitely generated as an
  $\mathcal{O}$-module.
  \item[(spl)] The field $K$ is a finite Galois extension of $F$ and
  the algebra $KA$ is split (i.e., for every simple $KA$-module $V$, $\mathrm{End}_{KA}(V) \simeq K$) semisimple.
\end{description}
\end{ypoth}
 
\begin{definition}\label{symmetric algebra}
We say that a linear map $t:A \rightarrow \mathcal{O}$ is a
symmetrizing form on $A$ or that $A$ is a symmetric algebra if 
\begin{enumerate}[(a)]
\item $t$ is a trace function, i.e., $t(ab)=t(ba)$ for all $a,b \in
A$,
\item the
morphism
$$\hat{t}:A \rightarrow \mathrm{Hom}_\mathcal{O}(A,\mathcal{O}),\,\,
  a \mapsto (x \mapsto \hat{t}(a)(x):=t(ax))$$
is an isomorphism of $A$-modules-$A$.
\end{enumerate}
\end{definition}

\begin{px}\label{symmetrizing form of the group algebra}
\small{\emph{In the case where $\mathcal{O}=\mathbb{Z}$ and
$A=\mathbb{Z}[G]$
 ($G$ a finite group), we can define the following symmetrizing form
 (``canonical'')
 on $A$
$$t:\mathbb{Z}[G] \rightarrow \mathbb{Z}, \,\,\, \sum_{g \in G}a_g g \mapsto a_1,$$
where $a_g \in \mathbb{Z}$ for all $g \in G$.}}
\end{px}

From now on, let us suppose that $A$ is a symmetric algebra with symmetrizing form $t$. 
By \cite{Ge}, we have the following results.

\begin{thedef}\label{schur elements and idempotents}\
\begin{enumerate}
  \item We have
  $$t=\sum_{\chi \in \mathrm{Irr}(KA)}\frac{1}{s_\chi}\chi,$$
  where $s_\chi$ is the Schur element of $\chi$
  with respect to $t$. 
  \item For all $\chi \in \mathrm{Irr}(KA)$, the central primitive
  idempotent associated to $\chi$ is
  $$e_\chi=\frac{1}{s_\chi}\sum_{i \in I}
  \chi(e_i')e_i,$$
  where $(e_i)_{i \in I}$ is a basis of $A$ over
$\mathcal{O}$ and $(e_i')_{i \in I}$ is the dual basis with respect
to $t$ (i.e., $t(e_ie_j')=\delta_{ij}$).
\end{enumerate}
\end{thedef}

\begin{corollary}\label{what we are searching}
The blocks of $A$ are the non-empty subsets $B$ of
$\mathrm{Irr}(KA)$ minimal with respect to the property
$$\sum_{\chi \in B}\frac{1}{s_\chi}\chi(a) \in \mathcal{O}, \textrm{ for all } a \in A.$$
\end{corollary}

Let us suppose now that $\mathcal{O}$ is a discrete valuation ring with unique prime ideal $\mathfrak{p}$ and that $K$ is the field of fractions of $\mathcal{O}$. Then the following result gives a criterion for a block to be a singleton.

\begin{proposition}\label{Malle-Rouquier}
Let $\chi \in \mathrm{Irr}(KA)$. The character $\chi$ is a block of $A$ by itself if
and only if $s_\chi \notin \mathfrak{p}.$
\end{proposition}
\begin{apod}{If $s_\chi \notin \mathfrak{p}$, then $1/s_\chi \in \mathcal{O}$ and
Corollary $\ref{what we are searching}$ implies that the character $\chi$ is a block of $A$ by itself. The inverse is a consequence of a theorem by Geck and Rouquier (cf. \cite{GeRo}, Proposition 4.4).}
\end{apod}

\subsection {Twisted symmetric algebras of finite groups}

Let $A$ be an $\mathcal{O}$-algebra such that the assumptions
$\ref{properties of the ring}$ are satisfied with a symmetrizing
form $t$. Let $\bar{A}$ be a subalgebra of $A$ free and of finite
rank as $\mathcal{O}$-module.

\begin{definition}\label{symmetric subalgebra}
We say that $\bar{A}$ is
a symmetric subalgebra of $A$, if it satisfies the following two
conditions:
\begin{enumerate}[(1)]
  \item $\bar{A}$ is free (of finite rank) as an $\mathcal{O}$-module and the
  restriction $\mathrm{Res}_{\bar{A}}^A(t)$ of the form $t$ to $\bar{A}$ is a symmetrizing form
  on $\bar{A}$,
  \item $A$ is free (of finite rank) as an $\bar{A}$-module for the action
  of left multiplication by the elements of $\bar{A}$.
\end{enumerate}
\end{definition}
We denote by
$$\mathrm{Ind}_{\bar{A}}^A: _{\bar{A}}\mathrm{\textbf{mod}} \rightarrow _A\mathrm{\textbf{mod}}
\,\textrm{ and }\, \mathrm{Res}_{\bar{A}}^A: _A\mathrm{\textbf{mod}}
\rightarrow _{\bar{A}}\mathrm{\textbf{mod}} $$ the functors defined
as usual by
$$\mathrm{Ind}_{\bar{A}}^A:=A \otimes_{\bar{A}}- \textrm{ where $A$ is viewed as an $A$-module-$\bar{A}$}$$
and
$$\mathrm{Res}_{\bar{A}}^A:=A \otimes_A - \textrm{ where $A$ is viewed as an $\bar{A}$-module-$A$.}$$ $ $

In the next sections, we will work on the Hecke algebras of complex
reflection groups, which are symmetric.
Sometimes the Hecke algebra of a group $W$ appears as a symmetric
subalgebra of the Hecke algebra of another group $W'$, which
contains $W$. Since we will be mostly interested in the
determination of the blocks of these algebras, it would be helpful,
if we could obtain the blocks of the former from the blocks of the
latter. This is possible with the use of a generalization of some
classic results, known as ``Clifford theory'' (see, for example,
\cite{Da}), to the twisted symmetric algebras of finite groups and
more precisely of finite cyclic groups. 

\begin{definition}\label{symmetric algebra of a finite group}
We say that a symmetric $\mathcal{O}$-algebra $(A,t)$ is the twisted
symmetric algebra of a finite group $G$ over the subalgebra
$\bar{A}$, if the following conditions are satisfied:
\begin{itemize}
  \item $\bar{A}$ is a symmetric subalgebra of $A$.
  \item There exists a family $\{A_g \,|\, g \in G\}$ of
  $\mathcal{O}$-submodules of $A$ such that
  \begin{description}
    \item[(a)] $A= \bigoplus_{g \in G}A_g$,
    \item[(b)] $A_gA_h=A_{gh}$ for all $g,h \in G$,
    \item[(c)] $A_1=\bar{A}$,
    \item[(d)] $t(A_g)=0$ for all $g \in G, g \neq  1$,
    \item[(e)] $A_g \cap A^\times \neq \emptyset$ for all $g \in G$ (where
    $A^\times$ is the set of units of $A$).
     \end{description}
    In particular, if $a_g \in A_g \cap A^\times$, then we have $A_g=a_g\bar{A}=\bar{A}a_g$.
\end{itemize}
\end{definition}

\subsection*{\normalsize Action of $G$ on $Z\bar{A}$}

From now on, we assume that $(A,t)$ is the twisted
symmetric algebra of a finite group $G$ over
$\bar{A}$ and that $K$ is an extension of $F$ such that
the algebras $KA$, $K\bar{A}$ and $KG$ are split semisimple.

\begin{thedef}\label{action of G on ZA}
Let $\bar{a} \in Z\bar{A}$ and $g \in G$. There exists a unique
element $g(\bar{a})$ of $\bar{A}$ satisfying
$$g(\bar{a})\mathfrak{g}=\mathfrak{g}\bar{a} \textrm{ for all } \mathfrak{g} \in
A_g.$$ If $a_g \in A^\times$ such that
$A_g=a_g\bar{A}$, then $$g(\bar{a})=a_g\bar{a}{a_g}^{-1}.$$
The map $\bar{a} \mapsto g(\bar{a})$ defines an action of $G$ as
ring automorphism of $Z\bar{A}$.
\end{thedef}

\subsection*{\normalsize Induction and restriction of $KA$-modules and
$K\bar{A}$-modules}

For all $\bar{\chi} \in \mathrm{Irr}(K\bar{A})$, we denote by
$\bar{e}(\bar{\chi})$ the block-idempotent of $K\bar{A}$ associated to
$\bar{\chi}$. If $g \in G$, then
$g(\bar{e}(\bar{\chi}))$ is also a block of $K\bar{A}$. 
Since
$K\bar{A}$ is split semisimple, it must be associated to an
irreducible character $g(\bar{\chi})$ of $K\bar{A}$. Thus, we can
define an action of $G$ on $\mathrm{Irr}(K\bar{A})$ such that for
all $g \in G$, $\bar{e}(g(\bar{\chi}))=g(\bar{e}(\bar{\chi}))$. We
denote by $G_{\bar{\chi}}$ the stabilizer of the character $\bar{\chi}$ in $G$ and by
$\bar{\Omega}$ the orbit of $\bar{\chi}$ under the action of $G$. We have
$|\bar{\Omega}|=|G|/|G_{\bar{\chi}}|$. We define
$$\bar{e}(\bar{\Omega}):=\sum_{g \in
G/G_{\bar{\chi}}}\bar{e}(g(\bar{\chi}))= \sum_{g \in
G/G_{\bar{\chi}}}g(\bar{e}(\bar{\chi})).$$ $ $\\
\emph{Case where $G$ is cyclic}\\
\\
Since the group $G$ is abelian, the set $\mathrm{Irr}(KG)$ forms a
group, which we denote by $G^\vee$. The application $\psi \mapsto
\psi \cdot \xi$, where $\psi \in \mathrm{Irr}(KA)$ and $\xi \in
G^\vee$, defines an action of $G^\vee$ on $\mathrm{Irr}(KA)$.
Then we have the following result

\begin{proposition}\label{1.42}
If the group $G$ is cyclic, there exists a bijection
$$\begin{array}{ccc}
    \mathrm{Irr}(K\bar{A})/G & \tilde{\leftrightarrow} & \mathrm{Irr}(KA)/G^\vee \\
    \bar{\Omega} & \leftrightarrow & \Omega
  \end{array}$$
such that
$$\bar{e}(\bar{\Omega})=e(\Omega),\, |\bar{\Omega}||\Omega|=|G| \textrm{ and }
\left\{
  \begin{array}{ll}
    \forall \chi \in \Omega, & \mathrm{Res}_{K\bar{A}}^{KA}(\chi)=\sum_{\bar{\chi} \in \bar{\Omega}}\bar{\chi}\\
    \forall \bar{\chi} \in \bar{\Omega}, &\mathrm{Ind}_{K\bar{A}}^{KA}(\bar{\chi})=\sum_{\chi \in \Omega}\chi
  \end{array}
\right.
$$
Moreover, for all $\chi \in \Omega$ and $\bar{\chi} \in
\bar{\Omega}$, we have
$$s_\chi = |\Omega| s_{\bar{\chi}}.$$
\end{proposition}

\subsection*{\normalsize Blocks of $A$ and blocks of $\bar{A}$}

Let us denote by $\mathrm{Bl}(A)$ the set of blocks of $A$ and by
$\mathrm{Bl}(\bar{A})$ the set of blocks of $\bar{A}$. For $\bar{b}
\in \mathrm{Bl}(\bar{A})$, we set
$$\mathrm{Tr}(G,\bar{b}):=\sum_{g \in G/G_{\bar{b}}}g(\bar{b}).$$
The algebra $(Z\bar{A})^G$ is contained in both $Z\bar{A}$ and $ZA$
and the set of its blocks is
$$\mathrm{Bl}((Z\bar{A})^G)=\{\mathrm{Tr}(G,\bar{b}) \,|\, \bar{b} \in
\mathrm{Bl}(\bar{A})/G\}.$$ Moreover, $\mathrm{Tr}(G,\bar{b})$ is a
sum of blocks of $A$ and we define the subset
$\mathrm{Bl}(A,\bar{b})$ of $\mathrm{Bl}(A)$ as follows
$$\mathrm{Tr}(G,\bar{b}):=\sum_{b \in \mathrm{Bl}(A,\bar{b})}b.$$

\begin{lemma}\label{1.43}
Let $\bar{b}$ be a block of $\bar{A}$ and
$\bar{B}:=\mathrm{Irr}(K\bar{A}\bar{b})$. Then
\begin{enumerate}[(1)]
  \item For all $\bar{\chi} \in \bar{B}$, we have $G_{\bar{\chi}}
  \subseteq G_{\bar{b}}$.
  \item We have
  $$\mathrm{Tr}(G,\bar{b})=\sum_{\bar{\chi} \in \bar{B}/G}
  \mathrm{Tr}(G,\bar{e}(\bar{\chi}))=
  \sum_{\{\bar{\Omega}|\bar{\Omega} \cap \bar{B} \neq \emptyset \}}\bar{e}(\bar{\Omega}).$$
\end{enumerate}
\end{lemma}

Now let $G^\vee:=\mathrm{Hom}(G,K^\times)$. We suppose that $K=F$.
The multiplication of the characters of $KA$ by the characters of
$KG$ defines an action of the group $G^\vee$ on $\mathrm{Irr}(KA)$.
This action is induced by the operation of $G^\vee$ on the algebra
$A$, which is defined in the following way:
$$ \xi \cdot (\bar{a}a_g) := \xi(g)\bar{a}a_g \,\,\textrm{ for all }
\xi \in G^\vee, \bar{a} \in \bar{A}, g \in G.$$ In particular,
$G^\vee$ acts on the set of blocks of $A$. Let $b$ be a block of
$A$. Denote by $\xi \cdot b$ the product of $\xi$ and $b$ and by
$(G^\vee)_b$ the stabilizer of $b$ in $G^\vee$. We set
$$\mathrm{Tr}(G^\vee,b):=\sum_{\xi\in G^\vee/(G^\vee)_b}\xi\cdot b.$$
The set of blocks of the algebra $(ZA)^{G^\vee}$ is given by
$$\mathrm{Bl}((ZA)^{G^\vee})=\{\mathrm{Tr}(G^\vee,b) \,|\, b \in
\mathrm{Bl}(A)/G^\vee\}.$$
The following lemma is the analogue of Lemma $\ref{1.43}$.

\begin{lemma}\label{1.44}
Let $b$ be a block of $A$ and $B:=\mathrm{Irr}(KAb)$. Then
\begin{enumerate}[(1)]
  \item For all $\chi \in B$, we have $(G^\vee)_\chi \subseteq (G^\vee)_b$.
  \item We have
  $$\mathrm{Tr}(G^\vee,b)=\sum_{\chi\in B/G^\vee} \mathrm{Tr}(G^\vee,e(\chi))=
  \sum_{\{\Omega|\Omega\cap B \neq \emptyset \}}e(\Omega).$$
\end{enumerate}
\end{lemma}
\emph{Case where $G$ is cyclic}\\
\\
For every orbit $\mathcal{Y}$ of $G^\vee$ on $\mathrm{Bl}(A)$, we
denote by $b(\mathcal{Y})$ the block of $(ZA)^{G^\vee}$ defined by
$$b(\mathcal{Y}):=\sum_{b \in \mathcal{Y}}b.$$
For every orbit $\bar{\mathcal{Y}}$ of $G$ on
$\mathrm{Bl}(\bar{A})$, we denote by $\bar{b}(\bar{\mathcal{Y}})$
the block of $(Z\bar{A})^G$ defined by
$$\bar{b}(\bar{\mathcal{Y}}):=\sum_{\bar{b} \in \bar{\mathcal{Y}}}\bar{b}.$$
The following proposition results from Proposition $\ref{1.42}$ and
Lemmas $\ref{1.43}$ and $\ref{1.44}$.

\begin{proposition}\label{1.45}
If the group $G$ is cyclic, there exists a bijection
$$\begin{array}{ccc}
    \mathrm{Bl}(\bar{A})/G & \tilde{\leftrightarrow} & \mathrm{Bl}(A)/G^\vee \\
    \bar{\mathcal{Y}} & \leftrightarrow & \mathcal{Y}
  \end{array}$$
such that
$$\bar{b}(\bar{\mathcal{Y}})=b(\mathcal{Y}),$$
i.e.,
$$\mathrm{Tr}(G,\bar{b})=\mathrm{Tr}(G^\vee,b) \textrm{ for all }
\bar{b} \in \bar{\mathcal{Y}} \textrm{ and } b \in \mathcal{Y}.$$ In
particular, the algebras $(Z\bar{A})^G$ and $(ZA)^{G^\vee}$ have the
same blocks.
\end{proposition}

\begin{corollary}\label{clifford}
If the blocks of $A$ are stable by the action of $G^\vee$, then the
blocks of $A$ coincide with the blocks of $(Z\bar{A})^G$.
\end{corollary}

\section {Hecke algebras of complex reflection groups}

\subsection{Generic Hecke algebras}

Let $\mu_\infty$ be the group of all the roots of unity in
$\mathbb{C}$ and $K$ a number field contained in
$\mathbb{Q}(\mu_\infty)$. We denote by $\mu(K)$ the group of all the
roots of unity of $K$. For every integer $d>1$, we set
$\zeta_d:=\mathrm{exp}(2\pi i/d)$ and denote by $\mu_d$ the group of
all the $d$-th roots of unity. 

Let $V$ be a $K$-vector space of
finite dimension $r$. Let $W$ be a finite subgroup of $\mathrm{GL}(V)$ generated by
(pseudo-)reflections acting irreducibly on $V$. Let us denote by $\mathcal{A}$ the set of the
reflecting hyperplanes of $W$. We set $\mathcal{M} := \mathbb{C} \otimes V - 
\bigcup_{H \in \mathcal{A}} \mathbb{C} \otimes H$. For $x_0 \in \mathcal{M}$, 
let $P:=\Pi_1(\mathcal{M},x_0)$ and $B:=\Pi_1(\mathcal{M}/W,x_0)$. Then there exists a short exact 
sequence (cf. \cite{BMR}):
$$\{1\}\rightarrow P \rightarrow B \rightarrow W \rightarrow\{1\}.$$
We denote by $\tau$ the central element of $P$ defined by the loop
$$[0,1] \rightarrow \mathcal{M}, \,\,\,\,\,\,\, t \mapsto \mathrm{exp}(2\pi it)x_0.$$ 

For every orbit $\mathcal{C}$ of $W$ on $\mathcal{A}$, we denote by
$e_{\mathcal{C}}$ the common order of the subgroups $W_H$, where $H$
is any element of $\mathcal{C}$ and $W_H$ the subgroup formed by $\mathrm{id}_V$
and all the reflections fixing the hyperplane $H$.

We choose a set of indeterminates
$\textbf{u}=(u_{\mathcal{C},j})_{(\mathcal{C} \in
\mathcal{A}/W)(0\leq j \leq e_{\mathcal{C}}-1)}$ and we denote by
$\mathbb{Z}[\textbf{u},\textbf{u}^{-1}]$ the Laurent polynomial ring
in all the indeterminates $\textbf{u}$. We define the \emph{generic
Hecke algebra} $\mathcal{H}$ of $W$ to be the quotient of the group
algebra $\mathbb{Z}[\textbf{u},\textbf{u}^{-1}]B$ by the ideal
generated by the elements of the form
$$(\textbf{s}-u_{\mathcal{C},0})(\textbf{s}-u_{\mathcal{C},1}) \cdots (\textbf{s}-u_{\mathcal{C},e_{\mathcal{C}}-1}),$$
where $\mathcal{C}$ runs over the set $\mathcal{A}/W$ and
$\textbf{s}$ runs over the set of monodromy generators around the
images in $\mathcal{M}/W$ of the elements of the hyperplane
orbit $\mathcal{C}$.

We make some assumptions for the algebra $\mathcal{H}$. Note that
they have been verified for all but a finite number of irreducible
complex reflection groups (\cite{BMM2}, remarks before 1.17, $\S$ 2;
\cite{GIM}).

\begin{ypoth}\label{ypo}
The algebra $\mathcal{H}$ is a free
$\mathbb{Z}[\textbf{\emph{u}},\textbf{\emph{u}}^{-1}]$-module of
rank $|W|$. Moreover, there exists a linear form
$t:\mathcal{H}\rightarrow
\mathbb{Z}[\textbf{\emph{u}},\textbf{\emph{u}}^{-1}]$ with the
following properties:
\begin{enumerate}[(1)]
    \item $t$ is a symmetrizing form on $\mathcal{H}$.
    \item Via the specialization $u_{\mathcal{C},j} \mapsto
     \zeta_{e_\mathcal{C}}^j$, the form $t$ becomes the canonical
     symmetrizing form on the group algebra $\mathbb{Z}W$.
    \item If we denote by $\alpha \mapsto \alpha^*$ the automorphism of
     $\mathbb{Z}[\emph{\textbf{u}},\emph{\textbf{u}}^{-1}]$ consisting of the
     simultaneous inversion of the indeterminates, then for all $b \in B$, we
     have
          $$t(b^{-1})^*=\frac{t(b\tau)}{t(\tau)}.$$
\end{enumerate}
\end{ypoth}

We know that the form $t$ is unique (\cite{BMM2}, 2.1). From now on,
let us suppose that the assumptions $\ref{ypo}$ are satisfied. Then
we have the following result by G.Malle (\cite{Ma4}, 5.2).

\begin{theorem}\label{Semisimplicity Malle}
Let $\textbf{\emph{v}}=(v_{\mathcal{C},j})_{(\mathcal{C} \in
\mathcal{A}/W)(0\leq j \leq e_{\mathcal{C}}-1)}$ be a set of
$\sum_{\mathcal{C} \in \mathcal{A}/W}e_{\mathcal{C}}$ indeterminates
such that, for every $\mathcal{C},j$, we have
$v_{\mathcal{C},j}^{|\mu(K)|}=\zeta_{e_\mathcal{C}}^{-j}u_{\mathcal{C},j}$.
Then the $K(\textbf{\emph{v}})$-algebra
$K(\textbf{\emph{v}})\mathcal{H}$ is split semisimple.
\end{theorem}

By ``Tits' deformation theorem'' (cf., for example, \cite{BMM2}, 7.2), it follows
that the specialization $v_{\mathcal{C},j}\mapsto 1$ induces a
bijection $\chi \mapsto \chi_{\textbf{v}}$ from the
set $\mathrm{Irr}(K(\textbf{v})\mathcal{H})$ of absolutely
irreducible characters of $K(\textbf{v})\mathcal{H}$ to the set
$\mathrm{Irr}(W)$ of absolutely irreducible characters of $W$.

The following result concerning the form of the Schur elements associated
with the irreducible characters of $K(\textbf{v})\mathcal{H}$ is proved in \cite{Chlou}, Theorem 4.2.5, using case by case analysis.

\begin{theorem}\label{Schur element generic}
The Schur element $s_\chi(\textbf{\emph{v}})$ associated with the
character $\chi_{\textbf{\emph{v}}}$ of
$K(\textbf{\emph{v}})\mathcal{H}$ is an element of
$\mathbb{Z}_K[\textbf{\emph{v}},\textbf{\emph{v}}^{-1}]$ of the form:
$$s_\chi({\textbf{\emph{v}}})=\xi_\chi N_\chi \prod_{i \in I_\chi} \Psi_{\chi,i}(M_{\chi,i})^{n_{\chi,i}}$$
where
\begin{itemize}
    \item $\xi_\chi$ is an element of $\mathbb{Z}_K$,
    \item $N_\chi= \prod_{\mathcal{C},j} v_{\mathcal{C},j}^{b_{\mathcal{C},j}}$ is a monomial in $\mathbb{Z}_K[\textbf{\emph{v}},\textbf{\emph{v}}^{-1}]$
          such that $\sum_{j=0}^{e_\mathcal{C}-1}b_{\mathcal{C},j}=0$
          for all $\mathcal{C} \in \mathcal{A}/W$,
    \item $I_\chi$ is an index set,
    \item $(\Psi_{\chi,i})_{i \in I_\chi}$ is a family of $K$-cyclotomic polynomials in one variable
           (i.e., minimal polynomials of the roots of unity over $K$),
    \item $(M_{\chi,i})_{i \in I_\chi}$ is a family of monomials in $\mathbb{Z}_K[\textbf{\emph{v}},\textbf{\emph{v}}^{-1}]$
          and if $M_{\chi,i} = \prod_{\mathcal{C},j} v_{\mathcal{C},j}^{a_{\mathcal{C},j}}$,
          then $\textrm{\emph{gcd}}(a_{\mathcal{C},j})=1$
          and $\sum_{j=0}^{e_\mathcal{C}-1}a_{\mathcal{C},j}=0$
          for all $\mathcal{C} \in \mathcal{A}/W$,
    \item ($n_{\chi,i})_{i \in I_\chi}$ is a family of positive integers.
\end{itemize}
This factorization is unique in $K[\textbf{\emph{v}},\textbf{\emph{v}}^{-1}]$. Moreover, the monomials
 $(M_{\chi,i})_{i \in I_\chi}$ are unique up to inversion, whereas the coefficient $\xi_\chi$ is unique up to multiplication by a root of unity.
\end{theorem}

Let $A:=\mathbb{Z}_K[\textbf{v},\textbf{v}^{-1}]$ and $\mathfrak{p}$ be a prime ideal of $\mathbb{Z}_K$. 

\begin{definition}\label{p-essential monomial}
Let  $M = \prod_{\mathcal{C},j} v_{\mathcal{C},j}^{a_{\mathcal{C},j}}$ be a monomial in $A$
such that $\textrm{\emph{gcd}}(a_{\mathcal{C},j})=1$. We say that $M$ is $\mathfrak{p}$-essential
for a character $\chi \in \mathrm{Irr}(W)$, if there exists a $K$-cyclotomic polynomial $\Psi$ such that
\begin{itemize}
\item $\Psi(M)$ divides $s_\chi(\textbf{\emph{v}})$.
\item $\Psi(1)  \in  \mathfrak{p}$.
\end{itemize}
We say that $M$ is $\mathfrak{p}$-essential 
for $W$, if there exists a character $\chi \in \mathrm{Irr}(W)$ such that
$M$ is $\mathfrak{p}$-essential for $\chi$.
\end{definition}

The following proposition (\cite{Chlou}, Proposition 3.1.3) gives a characterization of $\mathfrak{p}$-essential monomials, which plays an essential role in the proof of Theorem $\ref{main theorem}$. 

\begin{proposition}\label{p-essential}
Let  $M = \prod_{\mathcal{C},j} v_{\mathcal{C},j}^{a_{\mathcal{C},j}}$ be a monomial in $A$
such that $\textrm{\emph{gcd}}(a_{\mathcal{C},j})=1$. We set  $\mathfrak{q}_M:=(M-1)A +\mathfrak{p}A$.
Then
\begin{enumerate}
\item The ideal $\mathfrak{q}_M$ is a prime ideal of $A$. 
\item $M$ is $\mathfrak{p}$-essential for $\chi \in \mathrm{Irr}(W)$ if and only if
$s_\chi(\textbf{\emph{v}})/\xi_\chi \in \mathfrak{q}_M$.
\end{enumerate}
\end{proposition}

If $M$ is a $\mathfrak{p}$-essential monomial for $W$, then Theorem $\ref{main theorem}$ establishes a relation between the blocks of the algebra
$A_{\mathfrak{q}_M}\mathcal{H}$ and the Rouquier blocks. The following results concerning the blocks of $A_{\mathfrak{q}_M}\mathcal{H}$ are proven in \cite{Chlou}, Propositions $3.2.3$ and $3.2.5$. 

\begin{proposition}\label{not essential for block}
Let $M = \prod_{\mathcal{C},j} v_{\mathcal{C},j}^{a_{\mathcal{C},j}}$ be a monomial in $A$
such that $\textrm{\emph{gcd}}(a_{\mathcal{C},j})=1$ and  $\mathfrak{q}_M:=(M-1)A +\mathfrak{p}A$.
Then
\begin{enumerate}
\item If two irreducible characters are in the same block of $A_{\mathfrak{p}A}\mathcal{H}$, then they are in the same block of $A_{\mathfrak{q}_M}\mathcal{H}$.
\item If $C$ is a block of $A_{\mathfrak{p}A}\mathcal{H}$ and $M$ is not $\mathfrak{p}$-essential for any irreducible character in $C$, then $C$ is a block of $A_{\mathfrak{q}_M}\mathcal{H}$.
\end{enumerate}
\end{proposition}

\subsection {Cyclotomic Hecke algebras}

Let $y$ be an indeterminate. We set $q:=y^{|\mu(K)|}.$

\begin{definition}\label{specialization}
A cyclotomic specialization of $\mathcal{H}$ is a
$\mathbb{Z}_K$-algebra morphism $\phi:
\mathbb{Z}_K[\textbf{\emph{v}},\textbf{\emph{v}}^{-1}]\rightarrow
\mathbb{Z}_K[y,y^{-1}]$ with the following properties:
\begin{itemize}
  \item $\phi: v_{\mathcal{C},j} \mapsto y^{n_{\mathcal{C},j}}$ where
  $n_{\mathcal{C},j} \in \mathbb{Z}$ for all $\mathcal{C}$ and $j$.
  \item For all $\mathcal{C} \in \mathcal{A}/W$, and if $z$ is another
  indeterminate, the element of $\mathbb{Z}_K[y,y^{-1},z]$ defined by
  $$\Gamma_\mathcal{C}(y,z):=\prod_{j=0}^{e_\mathcal{C}-1}(z-\zeta_{e_\mathcal{C}}^jy^{n_{\mathcal{C},j}})$$
  is invariant by the action of $\textrm{\emph{Gal}}(K(y)/K(q))$.
\end{itemize}
\end{definition}

If $\phi$ is a cyclotomic specialization of $\mathcal{H}$,
the corresponding \emph{cyclotomic Hecke algebra} is the
$\mathbb{Z}_K[y,y^{-1}]$-algebra, denoted by $\mathcal{H}_\phi$,
which is obtained as the specialization of the
$\mathbb{Z}_K[\textbf{v},\textbf{v}^{-1}]$-algebra $\mathcal{H}$ via
the morphism $\phi$. It also has a symmetrizing form $t_\phi$
defined as the specialization of the canonical form $t$.\\
\\
\begin{remark} \emph{Sometimes we describe the morphism $\phi$ by the
formula}
$$u_{\mathcal{C},j} \mapsto \zeta_{e_\mathcal{C}}^j q^{n_{\mathcal{C},j}}.$$
\end{remark}

The following result is proved in \cite{Chlou}, Proposition 4.3.4.

\begin{proposition}\label{cyclotomic split semisimple}
The algebra $K(y)\mathcal{H}_\phi$ is split semisimple. 
\end{proposition}

For $y=1$ this algebra specializes to the group
algebra $KW$ (the form $t_\phi$ becoming the canonical form on the
group algebra). Thus, by ``Tits' deformation theorem'', the
specialization $v_{\mathcal{C},j} \mapsto 1$ induces the following bijections:
$$\begin{array}{ccccc}
    \textrm{Irr}(K(\textbf{v})\mathcal{H}) & \leftrightarrow & \textrm{Irr}(K(y)\mathcal{H}_\phi) & \leftrightarrow & \textrm{Irr}(W) \\
    \chi_{\textbf{v}} & \mapsto & \chi_{\phi} & \mapsto & \chi.
  \end{array}$$

\subsection {Rouquier blocks of the cyclotomic Hecke algebras}

\begin{definition}\label{Rouquier ring}
We call Rouquier ring of $K$ and denote by $\mathcal{R}_K(y)$ the
$\mathbb{Z}_K$-subalgebra of $K(y)$
$$\mathcal{R}_K(y):=\mathbb{Z}_K[y,y^{-1},(y^n-1)^{-1}_{n\geq 1}]$$
\end{definition}

Let $\phi: v_{\mathcal{C},j} \mapsto y^{n_{\mathcal{C},j}}$ be a
cyclotomic specialization and $\mathcal{H}_\phi$ the corresponding
cyclotomic Hecke algebra. The \emph{Rouquier blocks} of
$\mathcal{H}_\phi$ are the blocks of the algebra
$\mathcal{R}_K(y)\mathcal{H}_\phi$.\\
\\
\begin{remark}\emph{ If we set $q:=y^{|\mu(K)|}$, then the corresponding cyclotomic Hecke algebra $\mathcal{H}_\phi$ can be considered either over the ring $\mathbb{Z}_{K}[y,y^{-1}]$ or over the ring $\mathbb{Z}_{K}[q,q^{-1}]$. We define the Rouquier blocks of  $\mathcal{H}_\phi$ to be the blocks of  $\mathcal{H}_\phi$ defined over the Rouquier ring $\mathcal{R}_{K}(y)$ in $K(y)$. However, in other texts, as, for example, in \cite{BK}, the Rouquier blocks are determined over the Rouquier ring $\mathcal{R}_{K}(q)$ in $K(q)$. Since $\mathcal{R}_{K}(y)$ is the integral closure of $\mathcal{R}_{K}(q)$ in $K(y)$, Proposition \cite{BK}, $1.12$ establishes a relation between the blocks of
$\mathcal{R}_{K}(y)\mathcal{H}_\phi$ and the blocks of $\mathcal{R}_{K}(q)\mathcal{H}_\phi$. Moreover, in the case where $\mathcal{H}$ is an Ariki-Koike algebra (see section $3.2$),  they coincide (cf. \cite{Chlou2}, Proposition 3.6).}
\end{remark}
$ $\\

Set $\mathcal{O}:=\mathcal{R}_{K}(y)$ and let $\mathfrak{p}$ be a prime ideal of $\mathbb{Z}_K$.
The ring $\mathcal{O}$ is a Dedekind ring (cf., for example, \cite{Chlou}, Proposition 4.4.2) and hence, its localization  $\mathcal{O}_{\mathfrak{p}\mathcal{O}}$ at the prime ideal generated by $\mathfrak{p}$ is a discrete valuation ring. Following \cite{Chlou2}, Proposition 2.14, we have:

\begin{proposition}\label{Rouquier blocks and central characters}
Two characters $\chi,\psi \in \emph{Irr}(W)$ are in the same Rouquier block of $\mathcal{H}_\phi$
if and only if there exists a finite sequence
$\chi_0,\chi_1,\ldots,\chi_n \in \emph{Irr}(W)$ and a finite
sequence $\mathfrak{p}_1,\ldots,\mathfrak{p}_n$ of  prime
ideals of $\mathbb{Z}_K$ such that
\begin{itemize}
  \item $\chi_0=\chi$ and $\chi_n=\psi$,
  \item for all $j$ $(1\leq j \leq n)$,\,\,
         the characters $\chi_{j-1}$ and $\chi_{j}$
         belong to the same block of
         $\mathcal{O}_{\mathfrak{p}_j\mathcal{O}}\mathcal{H}_\phi.$
 \end{itemize}                   
\end{proposition}

The above proposition implies that if we know the blocks of the algebra $\mathcal{O}_{\mathfrak{p}\mathcal{O}}\mathcal{H}_\phi$ for every  prime
ideal of $\mathbb{Z}_K$,
then we know the Rouquier blocks of  $\mathcal{H}_\phi$. In order to determine the former, we
can use the following theorem (\cite{Chlou}, Theorem 3.3.2):

\begin{theorem}\label{main theorem}
Let $A:=\mathbb{Z}_K[\textbf{\emph{v}},\textbf{\emph{v}}^{-1}]$ and $\mathfrak{p}$ be a prime
ideal of $\mathbb{Z}_K$. 
Let $M_1,\ldots, M_k$ be all the
$\mathfrak{p}$-essential monomials for $W$ such that $\phi(M_j)=1$
for all $j=1,\ldots,k$. Set $\mathfrak{q}_0:=\mathfrak{p}A$,
$\mathfrak{q}_j:=\mathfrak{p}A+(M_j-1)A$ for $j=1,\ldots,k$ and
$\mathcal{Q}:=\{\mathfrak{q}_0,\mathfrak{q}_1,\ldots,\mathfrak{q}_k\}$.
Two irreducible characters $\chi,\psi \in \textrm{\emph{Irr}}(W)$
are in the same block of $\mathcal{O}_{\mathfrak{p}\mathcal{O}}\mathcal{H}_\varphi$ if
and only if there exist a finite sequence
$\chi_0,\chi_1,\ldots,\chi_n \in \textrm{\emph{Irr}}(W)$ and a
finite sequence $\mathfrak{q}_{j_1},\ldots,\mathfrak{q}_{j_n} \in
\mathcal{Q}$ such that
\begin{itemize}
  \item $\chi_0=\chi$ and $\chi_n=\psi$,
  \item for all $i$ $(1\leq i \leq n)$,  the characters $\chi_{i-1}$ and $\chi_i$ are
  in the same block of $A_{\mathfrak{q}_{j_i}}\mathcal{H}$.
\end{itemize}
\end{theorem}

Let $\mathfrak{p}$ be a prime ideal of $\mathbb{Z}_K$ and  $\phi: v_{\mathcal{C},j} \mapsto y^{n_{\mathcal{C},j}}$ a cyclotomic specialization.
If $M=\prod_{\mathcal{C},j}v_{\mathcal{C},j}^{a_{\mathcal{C},j}}$
is a $\mathfrak{p}$-essential monomial for $W$, then
$$\phi(M)=1 \Leftrightarrow \sum_{\mathcal{C},j}a_{\mathcal{C},j}n_{\mathcal{C},j}=0.$$
Set $m:=\sum_{\mathcal{C}\in \mathcal{A}/W}e_\mathcal{C}$. The
hyperplane defined in $\mathbb{C}^m$ by the relation
$$\sum_{\mathcal{C},j}a_{\mathcal{C},j}t_{\mathcal{C},j}=0,$$ where
$(t_ {\mathcal{C},j})_{ \mathcal{C},j}$ is a set of $m$
indeterminates, is called \emph{$\mathfrak{p}$-essential hyperplane}
for $W$. A hyperplane in $\mathbb{C}^m$ is called \emph{essential}
for $W$, if it is $\mathfrak{p}$-essential for some prime ideal
$\mathfrak{p}$ of $\mathbb{Z}_K$ (Respectively, a monomial is called \emph{essential}
for $W$, if it is $\mathfrak{p}$-essential for some prime ideal
$\mathfrak{p}$ of $\mathbb{Z}_K$).

\begin{definition}\label{assoc with hyp}Let
$\phi: v_{\mathcal{C},j} \mapsto y^{n_{\mathcal{C},j}}$ be a cyclotomic specialization  such that  the integers
$n_{\mathcal{C},j}$ belong to only one essential hyperplane $H$ (resp. to no essential hyperplane). We say that $\phi$ is a cyclotomic specialization associated with the essential hyperplane $H$ (resp. with no essential hyperplane). We call Rouquier blocks associated with the hyperplane
$H$ (resp. with no essential hyperplane) and denote by  $\mathcal{B}^H$ (resp. $\mathcal{B}^\emptyset)$ the partition of $\mathrm{Irr}(W)$ into Rouquier blocks of $\mathcal{H}_\phi$.
\end{definition}

With the help of the above definition and thanks to Proposition $\ref{Rouquier blocks and central characters}$ and Theorem $\ref{main theorem}$, we obtain the following characterization for the Rouquier blocks of a cyclotomic Hecke algebra.

\begin{proposition}\label{explain AllBlocks}
Let  $\phi: v_{\mathcal{C},j} \mapsto y^{n_{\mathcal{C},j}}$ be a cyclotomic specialization. If the integers
 $n_{\mathcal{C},j}$ belong to no essential hyperplane, then 
the Rouquier blocks of the cyclotomic Hecke algebra $\mathcal{H}_\phi$ coincide with the partition  $\mathcal{B}^\emptyset$.
Otherwise, two irreducible characters $\chi, \psi \in \textrm{\emph{Irr}}(W)$ belong to the same Rouquier block of  $\mathcal{H}_\phi$ if
and only if there exist a finite sequence
$\chi_0,\chi_1,\ldots,\chi_n \in \textrm{\emph{Irr}}(W)$ and a
finite sequence $H_1,\ldots,H_n$ of essential hyperplanes that the $n_{\mathcal{C},j}$ belong to such that
\begin{itemize}
  \item $\chi_0=\chi$ and $\chi_n=\psi$,
  \item for all $i$ $(1\leq i \leq n)$,  the characters $\chi_{i-1}$ and $\chi_i$ belong to $\mathcal{B}^{H_i}$.
\end{itemize}
 \end{proposition}

 \subsection{Functions $a$ and $A$}

Following the notations in \cite{BMM2}, 6B, for every element $P(y)
\in \mathbb{C}(y)$, we call
\begin{itemize}
  \item \emph{valuation of $P(y)$ at $y$} and denote by $\mathrm{val}_y(P)$ the order of $P(y)$
  at 0 (we have $\mathrm{val}_y(P)<0$ if 0 is a pole of $P(y)$ and $\mathrm{val}_y(P)>0$ if 0 is a zero of $P(y)$),
  \item \emph{degree of $P(y)$ at $y$} and denote by $\mathrm{deg}_y(P)$ the opposite of the
  valuation of $P(1/y)$.
\end{itemize}
Moreover, if $q:=y^{|\mu(K)|}$, then
$$\mathrm{val}_q(P):=\frac{\mathrm{val}_y(P)}{|\mu(K)|} \textrm{ and }
\mathrm{deg}_q(P):=\frac{\mathrm{deg}_y(P)}{|\mu(K)|}.$$ 
For $\chi
\in \mathrm{Irr}(W)$, we define
$$a_{\chi_\phi}:=\mathrm{val}_q(s_{\chi_\phi}(y)) \,\textrm{ and }\,
A_{\chi_\phi}:=\mathrm{deg}_q(s_{\chi_\phi}(y)).$$ 
The following
result is proven in \cite{BK}, Proposition 2.9.

\begin{proposition}\label{aA}\
Let $\chi,\psi \in \mathrm{Irr}(W)$. If $\chi_\phi$ and
        $\psi_\phi$ belong to the same Rouquier block, then
        $$a_{\chi_\phi}+A_{\chi_\phi}=a_{\psi_\phi}+A_{\psi_\phi}.$$
\end{proposition}

The values of the functions $a$ and $A$ can be calculated from the generic Schur elements. In order to explain how, we need to introduce the following symbols:
  
\begin{definition}\label{symbols} Let  $n \in \mathbb{Z}$. We set
\begin{itemize}
  \item $n^+ := \left\{
                      \begin{array}{ll}
                        n, & \hbox{if $n > 0$,} \\
                        0, & \hbox{if $n \leq 0$.}
                      \end{array}
                    \right.$ and\,  
  $(y^n)^+ := n^+$.
   \item $n^- = \left\{
                      \begin{array}{ll}
                        n, & \hbox{if $n < 0$,} \\
                        0, & \hbox{if $n \geq 0$.}
                      \end{array}
                    \right.$ \,and\,  
  $(y^n)^- := n^-$.
\end{itemize}
\end{definition}

Now let us fix $\chi \in
\mathrm{Irr}(W)$. Following the notations of Theorem $\ref{Schur element generic}$, the generic Schur 
element $s_\chi(\textbf{v})$
associated to $\chi$ is an element of
$\mathbb{Z}_K[\textbf{v},\textbf{v}^{-1}]$ of the form:
$$s_\chi(\textbf{v})=\xi_\chi N_\chi \prod_{i \in I_\chi} \Psi_{\chi,i}(M_{\chi,i})^{n_{\chi,i}}.\,\,\,\,\,\,\,\,\,\,\,\,\,\,\,\,(\dag)$$

We fix the factorization $(\dag)$ for $s_\chi(\textbf{v}).$ The following result is used in \cite{ChDeg} in order to obtain that the functions $a$ and $A$ are constant on the Rouquier blocks of the cyclotomic Hecke algebras of the exceptional complex reflection groups. 

\begin{proposition}\label{Aa formula}
Let $\phi:v_{\mathcal{C},j} \mapsto y^{n_{\mathcal{C},j}}$ be a
cyclotomic specialization. Then
\begin{itemize}
  \item $\mathrm{val}_y(s_{\chi_\phi}(y)) =
 \phi(N_\chi)^+ +\phi(N_\chi)^-+
  \sum_{i \in I_\chi}
  n_{\chi,i}\mathrm{deg}(\Psi_{\chi,i})(\phi(M_{\chi,i}))^-$.
  \item $\mathrm{deg}_y(s_{\chi_\phi}(y)) =\phi(N_\chi)^++\phi(N_\chi)^-+
  \sum_{i \in I_\chi}
  n_{\chi,i}\mathrm{deg}(\Psi_{\chi,i})(\phi(M_{\chi,i}))^+$.
\end{itemize}
\end{proposition}

\section {Rouquier blocks of the cyclotomic Hecke algebras of $G(de,e,r)$, $r>2$}

In \cite{Kim}, Kim determined the Rouquier blocks for the cyclotomic Hecke algebras of $G(de,e,r)$, following the method used in \cite{BK} for $G(d,d,r)$: Clifford theory to obtain the blocks of $G(de,e,r)$ from the blocks of $G(de,1,r)$. However, due to the incorrect determination of the Rouquier blocks for $G(de,1,r)$ in \cite{BK}, we will proceed here to some modifications to the results and their proofs. 
Moreover, in the next section, we'll explain why we have to distinguish the case where $r=2$ (more precisely, where $r=2$ and $e$ is even).

\subsection {Combinatorics}

Let $\el=(\el_1,\el_2,\ldots,\el_h)$ be a partition, \ie a finite decreasing sequence of positive integers
$\el_1 \geq \el_2 \geq \cdots \geq \el_h \geq 1.$
The integer
$|\el|:=\el_1+\el_2+\cdots+\el_h$
is called \emph{the size of $\el$}. We also say that $\lambda$ \emph{is a partition of }
$|\el|$.
The integer $h$ is called \emph{the height of $\el$} and we set $h_\el:=h$. To each partition $\el$ we associate its \emph{$\beta$-number}, $\eb_\el=(\eb_1,\eb_2,\ldots,\eb_h)$, defined by
$$\eb_1:=h+\el_1-1,\eb_2:=h+\el_2-2,\ldots,\eb_h:=h+\el_h-h.$$

\subsection*{\normalsize Multipartitions}

From now on,  let $d$ be a positive integer. Let $\el=(\el^{(0)},\el^{(1)},\ldots,\el^{(d-1)})$ be a $d$-partition, \ie a family of $d$ partitions indexed by the set $\{0,1,\ldots,d-1\}$. We set 
$$h^{(a)}:=h_{\el^{(a)}}, \,\,\, \eb^{(a)}:=\eb_{\el^{(a)}}$$
and we have
$$ \el^{(a)}=(\el_1^{(a)},\el_2^{(a)},\ldots,\el_{h^{(a)}}^{(a)}).$$
The integer
$$|\el|:=\sum_{a=0}^{d-1}|\el^{(a)}|$$
is called \emph{the size of $\el$}. We also say that $\lambda$ \emph{is a $d$-partition of}
$|\el|$.

\subsection*{\normalsize  Ordinary symbols}

If $\eb=(\eb_1,\eb_2,\ldots,\eb_h)$ is a sequence of positive integers such that $\eb_1>\eb_2>\cdots>\eb_h$ and $m$ is a positive integer, then the $m$-``shifted'' of $\eb$ is the sequence of numbers defined by
$$\eb[m]=(\eb_1+m,\eb_2+m,\ldots,\eb_h+m,m-1,m-2,\ldots,1,0).$$

Let  $\el=(\el^{(0)},\el^{(1)},\ldots,\el^{(d-1)})$ be a $d$-partition.  We call \emph{$d$-height of $\el$} the family $(h^{(0)},h^{(1)},\ldots,h^{(d-1)})$ and we define the 
 \emph{height of $\el$} to be the integer
$$h_\el:=\mathrm{max}\,\{h^{(a)} \,|\, (0 \leq a \leq d-1)\}.$$

\begin{definition}\label{ordinary standard symbol}
The ordinary standard symbol of $\el$ is the family of numbers defined by
$B_\el=(B^{(0)},B^{(1)},\ldots,B^{(d-1)}),$
where, for all $a$ $(0 \leq a \leq d-1)$, we have
$$B^{(a)}:=\eb^{(a)}[h_\el-h^{(a)}].$$
\end{definition}

The \emph{ordinary content} of a $d$-partition of ordinary standard symbol $B_\el$ is the multiset
$$\mathrm{Cont}_\el = B^{(0)} \cup B^{(1)} \cup \cdots \cup B^{(d-1)}.$$

\subsection*{\normalsize Charged symbols}

Assume that we have a given ``weight system'', \ie a family of integers
$$m:=(m^{(0)},m^{(1)},\ldots,m^{(d-1)}).$$

Let  $\el=(\el^{(0)},\el^{(1)},\ldots,\el^{(d-1)})$ be a $d$-partition.  We call \emph{$(d,m)$-charged height of $\el$} the family $(hc^{(0)},hc^{(1)},\ldots,hc^{(d-1)})$, where
$$hc^{(0)}:=h^{(0)}-m^{(0)},hc^{(1)}:=h^{(1)}-m^{(1)},\ldots,hc^{(d-1)}:=h^{(d-1)}-m^{(d-1)}.$$
We define the 
 \emph{$m$-charged height of $\el$} to be the integer
$$hc_\el:=\mathrm{max}\,\{hc^{(a)} \,|\, (0 \leq a \leq d-1)\}.$$

\begin{definition}\label{charged standard symbol}
The $m$-charged standard symbol of $\el$ is the family of numbers defined by
$Bc_\el=(Bc^{(0)},Bc^{(1)},\ldots,Bc^{(d-1)}),$
where, for all $a$ $(0 \leq a \leq d-1)$, we have
$$Bc^{(a)}:=\eb^{(a)}[hc_\el-hc^{(a)}].$$
\end{definition}
\begin{remark} \emph{The ordinary standard symbol corresponds to the weight system }
\begin{center}
$m^{(0)}=m^{(1)}=\cdots=m^{(d-1)}=0.$
\end{center}
\end{remark}

The \emph{$m$-charged content} of a $d$-partition of $m$-charged standard symbol $Bc_\el$ is the multiset
$$\mathrm{Contc}_\el = Bc^{(0)} \cup Bc^{(1)} \cup \cdots \cup Bc^{(d-1)}.$$

\subsection{Ariki-Koike algebras}

The group $G(d,1,r)$ is the group of all monomial $r \times r$ matrices with entries in $\mu_d$. It is isomorphic to the wreath product $\mu_d \wr
\mathfrak{S}_r$  and its field of definition is $K:=\mathbb{Q}(\zeta_d)$. Its irreducible characters are indexed by the $d$-partitions of $r$. If $\el$ is a $d$-partition of $r$, then we denote by $\chi_\el$ the corresponding irreducible character of $G(d,1,r)$.

The \emph{generic Ariki-Koike algebra} is the algebra $\mathcal{H}_{d,r}$ generated over the Laurent polynomial ring in $d+1$ indeterminates  
$$\mathbb{Z}[u_0,u_0^{-1},u_1,u_1^{-1},\ldots,u_{d-1},u_{d-1}^{-1},x,x^{-1}]$$
by the elements $\mathrm{\textbf{s}},\mathrm{\textbf{t}}_1,\mathrm{\textbf{t}}_2,\ldots,\mathrm{\textbf{t}}_{r-1}$ satisfying the relations
\begin{itemize}
\item $\mathrm{\textbf{s}}\mathrm{\textbf{t}}_1\mathrm{\textbf{s}}\mathrm{\textbf{t}}_1=\mathrm{\textbf{t}}_1\mathrm{\textbf{s}}\mathrm{\textbf{t}}_1\mathrm{\textbf{s}}$, $\mathrm{\textbf{s}}\mathrm{\textbf{t}}_j=\mathrm{\textbf{t}}_j\mathrm{\textbf{s}} \textrm{ for } j\neq 1$,
\item $\mathrm{\textbf{t}}_j\mathrm{\textbf{t}}_{j+1}\mathrm{\textbf{t}}_j=\mathrm{\textbf{t}}_{j+1}\mathrm{\textbf{t}}_j\mathrm{\textbf{t}}_{j+1}$,  $ \mathrm{\textbf{t}}_i\mathrm{\textbf{t}}_j=\mathrm{\textbf{t}}_j\mathrm{\textbf{t}}_i \textrm{ for } |i-j|>1$,
\item $(\mathrm{\textbf{s}}-u_0)(\mathrm{\textbf{s}}-u_1)\cdots(\mathrm{\textbf{s}}-u_{d-1})=(\mathrm{\textbf{t}}_j-x)(\mathrm{\textbf{t}}_j+1)=0$.
\end{itemize}

Let $$\phi : \left\{ 
\begin{array}{ll} 
u_j \mapsto \zeta_d^j q^{m_j}, &(0 \leq j <d),\\ 
x \mapsto q^n
\end{array} \right. 
$$ be a cyclotomic specialization for $\mathcal{H}_{d,r}$.
Thanks to Proposition $\ref{explain AllBlocks}$, in order to determine the Rouquier blocks of $(\mathcal{H}_{d,r})_\phi$ for any $\phi$, it suffices to determine the Rouquier blocks associated with no and each essential hyperplane for $G(d,1,r)$. Following \cite{Chlou2}, the essential hyperplanes for $G(d,1,r)$ are 
\begin{itemize}
\item $kN+M_s-M_t=0$, where $-r<k<r$  and $0 \leq s<t<d$ such that $\zeta_d^s-\zeta_d^t$ belongs to a prime ideal of 
 $\mathbb{Z}[\zeta_{d}]$,
\item $N=0$.
\end{itemize} 

We have proved that (cf. \cite{Chlou2}, Propositions 3.12, 3.15, 3.17)

\begin{theorem}\label{combination}\
\begin{enumerate}
\item The Rouquier blocks associated with no essential hyperplane are trivial.
\item Two irreducible characters $\chi_\el$ and $\chi_\mu$ belong to the same Rouquier block associated with the essential hyperplane $kN+M_s-M_t=0$  if and only if the following two conditions are satisfied:
\begin{itemize}
\item We have $\el^{(a)}=\mu^{(a)}$ for all $a \notin \{s,t\}.$
\item If $\el^{st}:=(\el^{(s)},\el^{(t)})$ and $\mu^{st}:=(\mu^{(s)},\mu^{(t)})$, then
$\mathrm{Contc}_{\el^{st}}= \mathrm{Contc}_{\mu^{st}}$ with respect to the weight system $(0,k)$.
\end{itemize}
\item Two irreducible characters $\chi_\el$ and $\chi_\mu$ belong to the same Rouquier block associated with the essential hyperplane $N=0$  if and only if  $|\lambda^{(a)}|=|\mu^{(a)}|$ for all $a=0,1,\ldots,d-1$. 
\end{enumerate}
\end{theorem}

Following Proposition $\ref{explain AllBlocks}$, the above theorem gives us an algorithm for the determination of the Rouquier blocks of any cyclotomic Ariki-Koike algebra (cf. \cite{Chlou2}, Theorem 3.18).

\subsection {Rouquier blocks for $G(de,e,r)$, $r>2$}

The group $G(de,e,r)$ is the group of all  $r \times r$ monomial matrices with entries in $\mu_{de}$  such that the product of all non-zero entries lies in $\mu_d$.

Following Ariki \cite{Ariki}, we define the 
 Hecke algebra of $G(de,e,r)$, $r >2$, to be the algebra $\mathcal{H}_{de,e,r}$ generated over the Laurent polynomial ring in $d+1$ indeterminates  
$$\mathbb{Z}[x_0,x_0^{-1},x_1,x_1^{-1},\ldots,x_{d-1},x_{d-1}^{-1},z,z^{-1}]$$
by the elements $a_0,a_1,\ldots,a_r$ satisfying the relations
\begin{itemize}
\item $(a_0-x_0)(a_0-x_1)\cdots(a_0-x_{d-1})=(a_j-z)(a_j+1)=0$ for $j=1,\ldots,r$,
\item $a_1a_3a_1=a_3a_1a_3$, $a_ja_{j+1}a_j=a_{j+1}a_ja_{j+1}$ for $j=2,\ldots,r-1$,
\item $a_1a_2a_3a_1a_2a_3=a_3a_1a_2a_3a_1a_2$,
\item $a_1a_j=a_ja_1$ for $j=4,\ldots,r$,
\item $a_i a_j=a_j a_i$  for $2 \leq i <j \leq r$ with $j-i>1$,
\item $a_0a_1a_2=(z^{-1}a_1a_2)^{2-e}a_2a_0a_1+(z-1)\sum_{k=1}^{e-2}(z^{-1}a_1a_2)^{1-k}a_0a_1=a_1a_2a_0$,
\item $a_0a_j=a_ja_0$ for $j=3,\ldots,r$.
\end{itemize}
Let $$\vartheta : \left\{ 
\begin{array}{ll} 
x_j \mapsto \zeta_d^j q^{m_j} &(0 \leq j <d),\\ 
y \mapsto q^n
\end{array} \right. 
$$
be a cyclotomic specialization for  $\mathcal{H}_{de,e,r}$. In order to determine the Rouquier blocks of 
$(\mathcal{H}_{de,e,r})_\vartheta$, we might as well consider the cyclotomic specialization
$$\phi : \left\{ 
\begin{array}{ll} 
x_j \mapsto \zeta_d^j q^{em_j} &(0 \leq j <d),\\ 
y \mapsto q^{en}.
\end{array} \right. 
$$
Since the integers $\{(m_j)_{0 \leq j <d},n\}$ and $\{(em_j)_{0 \leq j <d},en\}$ belong to the same essential hyperplanes for $G(de,e,r)$, Proposition $\ref{explain AllBlocks}$ implies that the Rouquier blocks of $(\mathcal{H}_{de,e,r})_\vartheta$ coincide with the Rouquier blocks of $(\mathcal{H}_{de,e,r})_\phi$.\\

We now consider the generic Ariki-Koike algebra $\mathcal{H}_{de,r}$ generated over the ring
$$\mathbb{Z}[u_0,u_0^{-1},u_1,u_1^{-1},\ldots,u_{de-1},u_{de-1}^{-1},x,x^{-1}]$$
by the elements $\mathrm{\textbf{s}},\mathrm{\textbf{t}}_1,\mathrm{\textbf{t}}_2,\ldots,\mathrm{\textbf{t}}_{r-1}$ satisfying the relations described in the definition of section $3.2$. Let 
$$\phi' : \left\{ 
\begin{array}{ll} 
u_j \mapsto \zeta_{de}^{j} q^{n_j} &(0 \leq j <de, n_j := m_{j \,\mathrm{mod}\,d}), \\ 
x \mapsto q^{en}
\end{array} \right. 
$$
be the ``corresponding'' cyclotomic specialization for  $\mathcal{H}_{de,r}$, \ie the specialization with respect to the weight system $$(m_0,m_1,\ldots,m_{d-1},m_0,m_1,\ldots,m_{d-1},\ldots,m_0,m_1,\ldots,m_{d-1}).$$
Set $\mathcal{H}:=(\mathcal{H}_{de,r})_{\phi'}$ and let $\bar{\mathcal{H}}$ be the subalgebra of $\mathcal{H}$ generated by
$$\textbf{s}^e, \tilde{\textbf{t}}_1:=\textbf{s}^{-1}\textbf{t}_1\textbf{s}, 
\textbf{t}_1, \textbf{t}_2, \ldots, \textbf{t}_{r-1}.$$
We have
$$\prod_{j=0}^{d-1}(\textbf{s}^e-\zeta_d^jq^{em_j})=(\tilde{\textbf{t}}_1-q^{en})(\tilde{\textbf{t}}_1+1)=
(\textbf{t}_i-q^{en})(\textbf{t}_i+1)=0 \textrm{ for } i=1,\ldots, r-1.$$
Then, by \cite{Ariki}, Proposition 1.16, we know that the algebra $(\mathcal{H}_{de,e,r})_\phi$ is isomorphic to the algebra $\bar{\mathcal{H}}$ via the morphism
$$a_0 \mapsto  \textbf{s}^e, a_1 \mapsto \tilde{\textbf{t}}_1, a_j \mapsto \textbf{t}_{j-1}\,\, (2 \leq j \leq r).$$

The following result is due to Kim (\cite{Kim}, Proposition 3.1).

\begin{proposition}\label{KIm}
The algebra $\mathcal{H}$ is a free $\bar{\mathcal{H}}$-module of rank $e$ with basis $\{1,\textbf{\emph{s}},\ldots,\textbf{\emph{s}}^{e-1}\}$, i.e., $$\mathcal{H}=\bar{\mathcal{H}} \oplus \textbf{\emph{s}} \bar{\mathcal{H}} \oplus \cdots \oplus \textbf{\emph{s}}^{e-1}\bar{\mathcal{H}}.$$
\end{proposition} 

By \cite{BMM2}, Proposition 1.18, the algebra $\mathcal{H}$ is symmetric and $\bar{\mathcal{H}}$ is a symmetric subalgebra of $\mathcal{H}$. In particular, following Definition $\ref{symmetric algebra of a finite group}$, $\mathcal{H}$ is the twisted symmetric algebra of the cyclic group of order $e$ over $\bar{\mathcal{H}}$ (since $\textbf{s}$ is a unit in $\mathcal{H}$). Therefore, we can apply Proposition $\ref{1.45}$ and obtain (using the notations of section $1.3$):

\begin{proposition}\label{first step}
If $G$ is the cyclic group of order $e$ and $K:=\mathbb{Q}(\zeta_{de})$, then the block-idempotents of $(Z\mathcal{R}_K(q)\bar{\mathcal{H}})^G$ coincide with the block-idempotents of $(Z\mathcal{R}_K(q)\mathcal{H})^{G^\vee}$, where $\mathcal{R}_K(q)$ is the Rouquier ring of $K$.
\end{proposition}

The action of the cyclic group $G^\vee$ of order $e$ on $\mathrm{Irr}(K(q)\mathcal{H})$ corresponds to the action generated by the cyclic permutation by $d$-packages on the $de$-partitions (cf., for example, \cite{Ma4}, \S4.A):
$$\begin{array}{rl}
\tau_d: &(\el^{(0)},\ldots,\el^{(d-1)},\el^{(d)},\ldots,\el^{(2d-1)},\ldots,\el^{(ed-d)},\ldots,\el^{(ed-1)})\\ \mapsto &(\el^{(ed-d)},\ldots,\el^{(ed-1)},\el^{(0)},\ldots,\el^{(d-1)},\ldots, \el^{(ed-2d)},\ldots,\el^{(ed-d-1)}).
\end{array}$$

More generally, the symmetric group $\mathfrak{S}_{de}$ acts naturally on the set of $de$-partitions of $r$ :  If $\tau \in \mathfrak{S}_{de}$ and
$\nu=(\nu^{(0)},\nu^{(1)},\ldots,\nu^{(de-1)})$ is a  $de$-partition of $r$, then 
$\tau(\nu):=(\nu^{(\tau(0))},\nu^{(\tau(1))},\ldots,\nu^{(\tau(de-1))})$.
The group $G^\vee$ is the cyclic subgroup of $\mathfrak{S}_{de}$ generated by the element
$$\tau_d=\prod_{j=0}^{d-1}\,\prod_{k=1}^{e-1}(j,j+kd).$$

Recall that $\mathcal{H}$ is the cyclotomic Ariki-Koike algebra of $G(de,1,r)$ corresponding to the weight system
 $$(m_0,m_1,\ldots,m_{d-1},m_0,m_1,\ldots,m_{d-1},\ldots,m_0,m_1,\ldots,m_{d-1}).$$
Following Proposition $\ref{explain AllBlocks}$, the Rouquier blocks of $\mathcal{H}$ are unions of the Rouquier blocks associated with the essential hyperplanes of the form
$$M_{j+kd}=M_{j+ld} \,\,(0 \leq j < d) (0 \leq k < l < e).$$
In order to show that the Rouquier blocks of $\mathcal{H}$ are stable under the action of $G^\vee$, it suffices to prove the following lemma: 

\begin{lemma}\label{exchange}
Let $\el$  be a $de$-partition of $r$, $j \in \{0,\ldots,d-1\}$ and $k \in \{1,\ldots,e-1\}$.
If $\mu=(j,j+kd)\el$, then
$\chi_\el$ and $\chi_\mu$ belong to the same Rouquier block of $\mathcal{H}$.
\end{lemma}
\begin{apod}{Suppose that $e=p_1^{a_1}p_2^{a_2}\cdots p_m^{a_m}$, where $p_i$ are prime numbers such that
$p_s \neq p_t$ for $s \neq t$. For $s \in \{1,2,\ldots,m\}$, we set
$c_s:= e/p_s^{a_s}.$ Then $\mathrm{gcd}(c_s)=1$ and by Bezout's theorem, there exist integers $(b_s)_{1 \leq s \leq m}$ such that $\sum_{s=1}^mb_sc_s=1$. Consequently,
$k=\sum_{s=1}^mkb_sc_s.$ We set $k_s:=kb_sc_s.$

For all  $s \in \{1,2,\ldots,m\}$, the element $1-\zeta_e^{c_s}$ belongs to the prime ideal of $\mathbb{Z}[\zeta_{de}]$ lying over the prime number $p_s$. So does $1-\zeta_e^{k_s}$. Now set
 $$l_0:=0 \textrm{ and } l_s:=\sum_{t=1}^{s}k_t \,\,(\mathrm{ mod }\,e).$$
We have that the element $\zeta_{de}^{j+l_{s-1}d}-\zeta_{de}^{j+l_{s}d}=\zeta_{de}^{j+l_{s-1}d}(1-\zeta_{e}^{k_{s}})$ belongs to the prime ideal of $\mathbb{Z}[\zeta_{de}]$ lying over the prime number $p_s$.
Therefore, the hyperplane $M_{j+l_{s-1}d}=M_{j+l_sd}$ is essential for $G(de,1,r)$.
Following the characterization of the Rouquier blocks associated with that hyperplane 
by Theorem $\ref{combination}$ and the fact that 
the ordinary content is stable under the action of a transposition, we obtain that the Rouquier blocks of $\mathcal{H}$ are stabilized by the action of $\sigma_s:=(j+l_{s-1}d,j+l_{s}d)$.
Set \begin{center} $\sigma:= \sigma_1 \circ \sigma_2 \circ \cdots \circ \sigma_{m-1} \circ \sigma_m\circ \sigma_{m-1} \circ \cdots \circ \sigma_2 \circ \sigma_1.$\end{center}
Then the characters $\chi_\el$ and $\chi_{\sigma(\el)}$ belong to the same Rouquier block of $\mathcal{H}$. It easy to check that $\sigma(\el)=\mu$.}
\end{apod}

Now the following result is immediate.

\begin{proposition}\label{second step}
If $\el$ is a $de$-partition of $r$, then the characters $\chi_\el$ and $\chi_{\tau_d(\el)}$
 belong to the same Rouquier block of $\mathcal{H}$. Therefore, the blocks of $\mathcal{R}_K(q)\mathcal{H}$ are stable under the action of $G^\vee$.
\end{proposition}

Thanks to the above result, Proposition $\ref{first step}$ now reads as follows:

\begin{corollary}\label{third step}
The block-idempotents of $(Z\mathcal{R}_K(q)\bar{\mathcal{H}})^G$ coincide with the block-idempotents of $\mathcal{R}_K(q)\mathcal{H}$.
\end{corollary}

Before we state our main result on the determination of the Rouquier blocks of $\mathcal{H}$, we will introduce the notion of ``$d$-stuttering $de$-partition'', following \cite{Kim}.

\begin{definition}\label{begayante}
Let $\el$ be a $de$-partition of $r$. We say that $\el$ is $d$-stuttering, if it's fixed by the action of 
$G^\vee$, i.e., if it's of the form
$$\el=(\el^{(0)},\ldots,\el^{(d-1)},\el^{(0)},\ldots,\el^{(d-1)},\ldots,\el^{(0)},\ldots,\el^{(d-1)}),$$
 where the first $d$ partitions are repeated $e$ times.
\end{definition} 

We are now ready to prove the main result:

\begin{theorem}\label{main result} Let $\el$ be a $de$-partition of $r$ and $\chi_\el$  the corresponding irreducible character of $G(de,1,r)$. We define $\mathrm{Irr}(K(q)\bar{\mathcal{H}})_\el$ to be the subset of $\mathrm{Irr}(K(q)\bar{\mathcal{H}})$ with the property:
$$\mathrm{Res}^{K(q)\mathcal{H}}_{K(q)\bar{\mathcal{H}}}\chi_\el=\sum_{\bar{\chi} \in \mathrm{Irr}(K(q)\bar{\mathcal{H}})_\el}\bar{\chi}.$$
Then
\begin{enumerate}
\item If $\el$ is $d$-stuttering and $\chi_\el$ is a block of $\mathcal{R}_K(q)\mathcal{H}$ by itself, then there are $e$ irreducible characters 
$(\bar{\chi})_{\bar{\chi} \in \mathrm{Irr}(K(q)\bar{\mathcal{H}})_\el}$. Each of these characters is a block of $\mathcal{R}_K(q)\bar{\mathcal{H}}$ by itself.
\item The other blocks of $\mathcal{R}_K(q)\mathcal{H}$ are in bijection with the blocks of  $\mathcal{R}_K(q)\bar{\mathcal{H}}$ via the map of Proposition $\ref{1.45}$, i.e., the corresponding block-idempotents of $\mathcal{R}_K(q)\mathcal{H}$ coincide with the remaining block-idempotents of $\mathcal{R}_K(q)\bar{\mathcal{H}}$.
\end{enumerate}
\end{theorem}
\begin{apod}{We will use here the notations of Propositions $\ref{1.42}$ and $\ref{1.45}$

If $\el$ is a $d$-stuttering partition, then it is the only element in its orbit $\Omega$ under the action of $G^\vee$. We have that $|\Omega||\bar{\Omega}|=|G|=e$, whence there exist $e$ elements in
$\bar{\Omega}= \mathrm{Irr}(K(q)\bar{\mathcal{H}})_\el$. If $\bar{\chi} \in \bar{\Omega}$, then its Schur element $s_{\bar{\chi}}$ is equal to the Schur element $s_\el$ of $\chi_\el$. If $\chi_\el$ is a block of $\mathcal{R}_K(q)\mathcal{H}$ by itself, then, by Propositions $\ref{Rouquier blocks and central characters}$ and
 $\ref{Malle-Rouquier}$, $s_\el$ is invertible in $\mathcal{R}_K(q)$ and so is $s_{\bar{\chi}}$. Thus, $\bar{\chi}$ is a block of $\mathcal{R}_K(q)\bar{\mathcal{H}}$ by itself. 

If $\el$ is not a $d$-stuttering partition and $b$ is the block containing  $\chi_\el$, then, in order to establish the desired bijection, we have to show that the block $\bar{b}$ of $\mathcal{R}_K(q)\bar{\mathcal{H}}$ which contains  a character in $\mathrm{Irr}(K(q)\bar{\mathcal{H}})_\el$ is fixed by the action of $G$, \ie that $\bar{b}=\mathrm{Tr}(G,\bar{b})$. Thanks to the  lemma that follows this theorem, for all prime divisor $p$ of $e$, there exists a
$de$-partition $\el(p)$ of $r$ such that $\chi_{\el(p)}$ belongs to $b$ and  the order of $G^\vee_{\chi_{\el(p)}}$ is not divisible by $p$. By Proposition $\ref{1.42}$, we know that for each $\bar{\chi} \in
\mathrm{Irr}(K(q)\bar{\mathcal{H}})_{\el(p)}$, we have $|G^\vee_{\chi_{\el(p)}}||G_{\bar{\chi}}|=e$. Thus, $|G_{\bar{\chi}}|$ is divisible by the largest power of $p$ dividing $e$.  Since $b=\mathrm{Tr}(G,\bar{b})$,  the elements of $\mathrm{Irr}(K(q)\bar{\mathcal{H}})_{\el(p)}$ belong to blocks of $\mathcal{R}_K(q)\bar{\mathcal{H}}$ conjugate of $\bar{b}$ by $G$, whose stabilizer is $G_{\bar{b}}$. By Lemma $\ref{1.43}$(1), we obtain that, for every prime number $p$,  $|G_{\bar{b}}|$ is divisible by the largest power of $p$ dividing $e$.
Thus, $G_{\bar{b}}=G$ and $\mathrm{Tr}(G,\bar{b})=\bar{b}$.

It remains to show that if $\el$ is a $d$-stuttering partition and $\chi_\el$ is not a block of $\mathcal{R}_K(q)\mathcal{H}$ by itself, then there exists a partition $\mu$ such that $\chi_\el$ and $\chi_{\mu}$ belong to the same block of $\mathcal{R}_K(q)\mathcal{H}$ and $\mu$ is not $d$-stuttering. Then the second case described above covers our needs.

If $\el$ is a $d$-stuttering partition, then the description of the Schur elements for $\mathcal{H}$
(cf., for example, \cite{Mat}, Corollary 6.5)  implies that  the essential hyperplanes of the form 
$$M_{j+kd}=M_{j+ld} \,\,(0 \leq j < d) (0 \leq k<l < e),$$
are not essential for $\chi_\el$. If now $\chi_\el$ is not a block of $\mathcal{R}_K(q)\mathcal{H}$ by itself, then, by Proposition $\ref{explain AllBlocks}$, there exists a $de$-partition $\mu \neq \el$ such that $\chi_\el$ and $\chi_\mu$ belong to the same Rouquier block associated with another essential hyperplane $H$ for $G(de,1,r)$ such that the integers $\{(n_j)_{0 \leq j < de}, en\}$ belong to $H$.

If $H$ is $N=0,$ then, by Theorem $\ref{combination}$, we have $|\el^{(a)}|=|\mu^{(a)}|$ for all $a=0,1,\ldots,de-1$. Since 
$\el \neq \mu$, there exists $s \in \{0,1,\ldots,de-1\}$ such that $ \el^{(s)} \neq \mu^{(s)}$. If  $\nu$ is the partition obtained from $\el$ by exchanging $\el^{(s)}$ and $\mu^{(s)}$, then $\chi_\el$ and $\chi_{\nu}$ belong to the same block of $\mathcal{R}_K(q)\mathcal{H}$ and $\nu$ is not $d$-stuttering.

If $H$ is of the form 
$kN+M_s-M_t=0, \textrm{ where } -r<k<r \textrm{ and } 0 \leq s<t<de$, then
$\el^{(a)}=\mu^{(a)}$  for all $a \neq s,t$. If $s \not\equiv t \,\mathrm{ mod }\,d$ or $e>2$, then
$\mu$ can not be $d$-stuttering. Suppose now that $s \equiv t \,\mathrm{ mod }\,d$ and $e=2$. 
As mentioned above, 
 the hyperplane $M_s=M_t$ is not essential for $\chi_\el$, whence $k \neq 0$. 
Since the integers $\{(n_j)_{0 \leq j < de}, en\}$ belong to $H$ and $n_s=n_t$, we must have $n=0$.
If $\mu$ is $d$-stuttering, then $\mu^{(s)}=\mu^{(t)}$ and we deduce that $|\mu^{(s)}|=|\mu^{(t)}|=|\el^{(t)}|=|\el^{(s)}|$.
Let $\nu$ be the $de$-partition obtained from $\el$ by replacing $\el^{(t)}$ with $\mu^{(t)}$. Then $\nu$ is not $d$-stuttering and the characters $\chi_\el$ and $\chi_{\nu}$ belong to the same Rouquier block associated with the essential hyperplane $N=0$. Since $n=0$,  Proposition $\ref{explain AllBlocks}$ implies that  $\chi_\el$ and $\chi_{\nu}$ belong to the same block of $\mathcal{R}_K(q)\mathcal{H}$.}
\end{apod}

\begin{lemma}\label{last lemma}
If  $\el$ is not a $d$-stuttering partition of $r$ and $p$ is a prime divisor of $e$, then there exists a
$de$-partition $\el(p)$ of $r$ such that $\chi_\el$ and $\chi_{\el(p)}$ belong to the same block of
$\mathcal{R}_K(q)\mathcal{H}$ and  the order of $G^\vee_{\chi_{\el(p)}}$ is not divisible by $p$.
\end{lemma}
\begin{apod}{If $\el=(\el^{(0)},\ldots,\el^{(d-1)},\el^{(d)},\ldots,\el^{(2d-1)},\ldots,\el^{(ed-d)},\ldots,\el^{(ed-1)}),$
then, for $i=0,1,\ldots,e-1$, we define the $d$-partition $\el_i$ as follows:
$$\el_i:=(\el^{(id)},\el^{(id+1)},\ldots,\el^{(id+d-1)}).$$
Then $\el=(\el_0,\el_1,\ldots,\el_{e-1}).$ Since $\el$ isn't $d$-stuttering, there exists $m \in \{ 0, 1,\ldots,e-1\}$ such that $\el_0 \neq \el_m$. We denote by $\el(p)$ the partition obtained from $\el$ by exchanging $\el_m$ and $\el_{e/p}$. Due to Lemma $\ref{exchange}$, the characters $\chi_\el$ and
$\chi_{\el(p)}$ belong to the same block of $\mathcal{R}_K(q)\mathcal{H}$.
Moreover, by construction, the $de$-partition $\el(p)$ isn't fixed by the generator of the unique subgroup of order $p$ of $G^\vee$, which proves that the order of its stabilizer is prime to $p$.}
\end{apod}

\subsection*{\normalsize Functions $a$ and $A$}

\begin{itemize}
\item The description of the Rouquier blocks of $\bar{\mathcal{H}}$ by Theorem $\ref{main result}$,
\item the relation between the Schur elements of $\bar{\mathcal{H}}$
and the Schur elements of $\mathcal{H}$ given by Proposition $\ref{1.42}$
\item and the invariance of the integers $a_\chi$ and $A_\chi$ on the Rouquier blocks of $\mathcal{H}$, resulting from propositions \cite{BK}, 3.18, and \cite{Chlou2}, 3.21  imply that
\end{itemize}

\begin{proposition}\label{fin}
The valuations $a_{\bar{\chi}}$ and the degrees $A_{\bar{\chi}}$ of the Schur elements are constant on the Rouquier blocks of $\bar{\mathcal{H}}$.
\end{proposition}

\section {Rouquier blocks of the cyclotomic Hecke algebras of $G(de,e,2)$}

If the integer $e$ is odd, then the Hecke algebra of the group $G(de,e,2)$ can be viewed as a symmetric subalgebra of a Hecke algebra of the group $G(de,1,2)$ and all the results of the previous section hold.\\

If $e$ is even, this can't be done, because there exist three orbits of reflecting hyperplanes under the action of the group. Following \cite{Ariki}, Proposition 1.16, Malle shows (cf.  \cite{Ma2}, Proposition $3.9$) that the Hecke algebra of the group 
$G(de,e,2)$ can be viewed as a symmetric subalgebra of a Hecke algebra of the group $G(de,2,2)$ and thus, we can apply Clifford theory in order to obtain the blocks of the former from the blocks of the latter. 

\subsection{Rouquier blocks for $G(2d,2,2)$}

Let $d\geq 1$. The group $G(2d,2,2)$ has $4d$ irreducible characters of degree $1$, 
$$\chi_{ijk} \,\,(0 \leq i,j\leq 1)\,(0 \leq k < d),$$
and $d^2-d$ irreducible characters of degree $2$,
$$\chi_{kl}^{1},\, \chi_{kl}^{2}\,\,(0 \leq k \neq l < d),$$
with $\chi_{kl}^{1,2}=\chi_{lk}^{1,2}$.\\

The generic Hecke algebra of the group $G(2d,2,2)$ is the algebra $\mathcal{H}_{d}$ generated over the Laurent  polynomial ring in $d+4$ indeterminates  
$$\mathbb{Z}[x_0,x_0^{-1},x_1,x_1^{-1},y_0,y_0^{-1},y_1,y_1^{-1},z_0,z_0^{-1},z_1,z_1^{-1},\ldots,z_{d-1},z_{d-1}^{-1}]$$
by the elements $s,t,u$ satisfying the relations
\begin{itemize}
\item $stu=tus=ust$,
\item $(s-x_0)(s-x_1)=(t-y_0)(t-y_1)=(u-z_0)(u-z_1)\cdots(u-z_{d-1})=0$.
\end{itemize}

The following theorem (\cite{Ma2}, Theorem $3.11$) gives a  description of the generic Schur elements for $G(2d,2,2)$.

\begin{theorem}\label{nice Schur}
Let us denote by $\Phi_1$ the first $\mathbb{Q}$-cyclotomic polynomial (i.e., $\Phi_1(q)=q-1$). The generic Schur elements for $\mathcal{H}_d$ are given by
\begin{center}
$\Phi_1(x_ix_{1-i}^{-1}) \cdot \Phi_1(y_jy_{1-j}^{-1})\cdot \prod_{l=0,\,l\neq k}^{d-1} (\Phi_1(z_kz_l^{-1})
\cdot \Phi_1(x_ix_{1-i}^{-1} y_jy_{1-j}^{-1}z_kz_l^{-1}))$
\end{center}
for the linear characters $\chi_{ijk}$, and
\begin{center}
$-2 \cdot \prod_{m=0,\,m\neq k,l}^{d-1}( \Phi_1(z_kz_m^{-1})\cdot \Phi_1(z_lz_m^{-1})) \cdot$ 
$\prod_{i=0}^1 (\Phi_1(X_iX_{1-i}^{-1}Y_iY_{1-i}^{-1}Z_kZ_l^{-1})\cdot
\Phi_1(X_iX_{1-i}^{-1}Y_{1-i}Y_{i}^{-1}Z_lZ_k^{-1})),$
\end{center}
with $X_i^2:=x_i$, $Y_j^2:=y_j$, $Z_k^2:=z_k$, for the characters $\chi_{kl}^{1,2}$ of degree $2$.
\end{theorem}

The field of definition of $G(2d,2,2)$ is $K:=\mathbb{Q}(\zeta_{2d})$. Following Theorem $\ref{Semisimplicity Malle}$, if we set
$$\mathcal{X}_i^{|\mu(K)|}:=(-1)^{-i}x_i \textrm{ for } i=0,1,\,\,\,  \mathcal{Y}_j^{|\mu(K)|}:=(-1)^{-j}y_j \textrm{ for } j=0,1$$
and
$$\mathcal{Z}_k^{|\mu(K)|}:=\zeta_d^{-k}z_j  \textrm{ for } k=0,1,\ldots,d-1,$$
then the algebra $K(\mathcal{X}_0,\mathcal{X}_1,\mathcal{Y}_0,\mathcal{Y}_1,\mathcal{Z}_0,\mathcal{Z}_1,\ldots,\mathcal{Z}_{d-1})\mathcal{H}_{d}$ is split semisimple.\\

Let $\mathfrak{I}$ be the prime ideal of $\mathbb{Z}[\zeta_{2d}]$ lying over $2$. The description of the generic Schur elements by Theorem $\ref{nice Schur}$ implies that the essential monomials for $G(2d,2,2)$ are
\begin{itemize}
\item $\mathcal{X}_0\mathcal{X}_1^{-1}$ ($\mathfrak{I}$-essential),
\item $\mathcal{Y}_0\mathcal{Y}_1^{-1}$ ($\mathfrak{I}$-essential),
\item $\mathcal{Z}_k\mathcal{Z}_l^{-1}$, where $0 \leq k<l<d$ such that $\zeta_d^k-\zeta_d^l$ belongs to a prime ideal $\mathfrak{p}$ of 
 $\mathbb{Z}[\zeta_{2d}]$ ($\mathfrak{p}$-essential),
\item $\mathcal{X}_i\mathcal{X}_{1-i}^{-1}\mathcal{Y}_j\mathcal{Y}_{1-j}^{-1}\mathcal{Z}_k\mathcal{Z}_l^{-1}$, where $0 \leq i,j \leq 1$ and $0 \leq k<l<d$ such that $\zeta_d^k-\zeta_d^l$ belongs to a prime ideal $\mathfrak{p}$ of $\mathbb{Z}[\zeta_{2d}]$ ($\mathfrak{p}$-essential).
\end{itemize}

Let $\phi$ be a cyclotomic specialization for $\mathcal{H}_d$, \ie a $\mathbb{Z}_K$-algebra morphism of the form
$$\phi: \mathcal{X}_i \mapsto y^{a_i}, \,\,
 \mathcal{Y}_j \mapsto y^{b_j}, \,\,
 \mathcal{Z}_k \mapsto y^{c_k}.$$ 
Set $q:=y^{|\mu(K)|}$. Then $\phi$ can be described as follows:
$$\phi: x_i \mapsto (-1)^iq^{a_i}, \,\,
 y_j \mapsto (-1)^jq^{b_j}, \,\,
 z_k \mapsto \zeta_d^kq^{c_k}.$$ 

Due to Proposition $\ref{cyclotomic split semisimple}$, "Tits' deformation theorem" implies that the specialization $y \mapsto 1$ induces a bijection
$$\begin{array}{ccc}
   \textrm{Irr}(K(y)(\mathcal{H}_d)_\phi) & \leftrightarrow & \textrm{Irr}(G(2d,2,2)) \\
  \chi_{\phi} & \mapsto & \chi.
  \end{array}$$
For $\chi \in \mathrm{Irr}(G(2d,2,2))$, let $s_{\chi_\phi}$  be the corresponding cyclotomic Schur element.
As in section $2.4$, we set
$$a_{\chi_\phi}:=\mathrm{val}_q(s_{\chi_\phi}(y))=\frac{\mathrm{val}_y(s_{\chi_\phi}(y))}{|\mu(K)|} \textrm{ and }
A_{\chi_\phi}:=\mathrm{deg}_q(s_{\chi_\phi}(y))=\frac{\mathrm{deg}_y(s_{\chi_\phi}(y))}{|\mu(K)|} .$$
Then, by Proposition $\ref{aA}$, we have that if two irreducible characters 
 $\chi_\phi$ and $\psi_\phi$ belong to the same Rouquier block of $(\mathcal{H}_d)_\phi$, then
        $$a_{\chi_\phi}+A_{\chi_\phi}=a_{\psi_\phi}+A_{\psi_\phi}.$$
Thanks to the formulas of Proposition $\ref{Aa formula}$, the following result derives immediately from the description of the generic Schur elements by Theorem $\ref{nice Schur}$.

\begin{proposition}\label{values of aA}
Let $\chi \in \mathrm{Irr}(G(2d,2,2))$. If $\chi$ is a linear character $\chi_{ijk}$, then
$$a_{\chi_\phi}+A_{\chi_\phi}=d(a_i-a_{1-i}+b_j-b_{1-j}+2c_k)-2\sum_{l=0}^{d-1}c_l.$$
If $\chi$ is a character $\chi_{kl}^{1,2}$ of degree $2$, then
$$a_{\chi_\phi}+A_{\chi_\phi}=d(c_k+c_l)-2 \sum_{m=0}^{d-1}c_m.$$
\end{proposition}

Following Proposition $\ref{explain AllBlocks}$, in order to determine the Rouquier blocks of the cyclotomic Hecke algebras of $G(2d,2,2)$, it suffices to determine the Rouquier blocks associated with its essential hyperplanes. 

\begin{theorem}\label{yes proof} For the group $G(2d,2,2)$, we have that
\begin{enumerate}[(1)]
\item The non-trivial Rouquier blocks associated with no essential hyperplane are 
$$\{\chi_{kl}^{1},\chi_{kl}^{2}\} \,\textrm{ for all } 0 \leq k<l<d.$$
\item The non-trivial Rouquier blocks associated with the $\mathfrak{I}$-essential hyperplane $A_0=A_1$ are
$$\{\chi_{0jk},\chi_{1jk}\}\, \textrm{ for all } 0 \leq j \leq 1 \textrm{ and } 0 \leq k<d,$$
$$\{\chi_{kl}^{1},\chi_{kl}^{2}\}\, \textrm{ for all } 0 \leq k<l<d.$$
\item The non-trivial Rouquier blocks associated with the $\mathfrak{I}$-essential hyperplane $B_0=B_1$ are
$$\{\chi_{i0k},\chi_{i1k}\} \,\textrm{ for all } 0 \leq i \leq 1 \textrm{ and } 0 \leq k<d,$$
$$\{\chi_{kl}^{1},\chi_{kl}^{2}\}\, \textrm{ for all } 0 \leq k<l<d.$$
\item The non-trivial Rouquier blocks associated with the $\mathfrak{p}$-essential hyperplane $C_k=C_l$ $(0 \leq k<l<d)$ are
$$ \{\chi_{ijk},\chi_{ijl}\}\, \textrm{ for all } 0 \leq i,j \leq 1,$$
$$\{\chi_{km}^{1},\chi_{km}^{2},\chi_{lm}^{1},\chi_{lm}^{2}\}\, \textrm{ for all } 0 \leq m <d \textrm{ with } m \notin \{k,l\},$$
$$\{\chi_{kl}^{1},\chi_{kl}^{2}\},$$
$$\{\chi_{rs}^{1},\chi_{rs}^{2}\}\, \textrm{ for all } 0 \leq r<s<d \textrm{ with } r,s \notin \{k,l\}.$$
\item The non-trivial Rouquier blocks associated with the $\mathfrak{p}$-essential hyperplane $A_i-A_{1-i}+B_j-B_{1-j}+C_k-C_l=0$ $(0 \leq i,j\leq 1)$ $(0 \leq k<l<d)$ are
$$\{\chi_{ijk},\chi_{1-i,1-j,l}, \chi_{kl}^{1},\chi_{kl}^{2} \}, $$
$$\{\chi_{rs}^{1},\chi_{rs}^{2}\}\, \textrm{ for all } 0 \leq r<s<d \textrm{ with } (r,s)  \neq (k,l).$$
\end{enumerate}
\end{theorem}
\begin{apod}{Following Definition $\ref{assoc with hyp}$, in each case, we need to determine the Rouquier blocks of a cyclotomic Hecke algebra obtained via a specialization associated with the corresponding essential hyperplane. We recall that, due to Proposition $\ref{Malle-Rouquier}$, if a hyperplane is essential for an irreducible character $\chi$, then $\chi$ isn't a Rouquier block by itself.
Moreover, the first part of Proposition $\ref{not essential for block}$ implies that the Rouquier blocks associated with an essential hyperplane are unions of the Rouquier blocks associated with no essential hyperplane.
\begin{enumerate}[(1)]
\item Let $\phi$ be any cyclotomic specialization associated with no essential hyperplane. Due to Proposition $\ref{Malle-Rouquier}$, each linear character is a Rouquier block by itself, whereas any character of degree $2$ isn't. Now, by Proposition $\ref{aA}$, we have that if two irreducible characters 
 $\chi_\phi$ and $\psi_\phi$ belong to the same Rouquier block of $(\mathcal{H}_d)_\phi$, then
 $a_{\chi_\phi}+A_{\chi_\phi}=a_{\psi_\phi}+A_{\psi_\phi}.$ The formulas of Proposition $\ref{values of aA}$ imply that the character $\chi_{kl}^1$ $(0 \leq k<l<d)$ can  be in the same block only with the character $\chi_{kl}^2$.
\item Let $\phi$ be any cyclotomic specialization associated with the $\mathfrak{I}$-essential hyperplane
$A_0=A_1$. Since this isn't an essential hyperplane for the characters of degree $2$,
Proposition $\ref{not essential for block}$ implies that $\{\chi_{kl}^{1},\chi_{kl}^{2}\}$ is a Rouquier block of $(\mathcal{H}_d)_\phi$ for all $0 \leq k<l<d$. Now, the hyperplane $A_0=A_1$ is $\mathfrak{I}$-essential for all characters of degree $1$ and thus, by Proposition $\ref{Malle-Rouquier}$, the linear characters don't form blocks by themselves. Due to Proposition $\ref{aA}$, the  formulas of Proposition $\ref{values of aA}$ imply that the character $\chi_{0jk}$ $(0 \leq j\leq 1,0 \leq k<d)$ can be in the same block only with the character $\chi_{1jk}$.
\item For the  $\mathfrak{I}$-essential hyperplane $B_0=B_1$, we use the same method as in the previous case.
\item Let $\phi$ be a cyclotomic specialization associated with the $\mathfrak{p}$-essential hyperplane
$C_k=C_l$, where $0 \leq k<l<d$. Since the Rouquier blocks associated with an essential hyperplane are unions of the Rouquier blocks associated with no essential hyperplane, the characters 
 $\chi_{rs}^{1}$ and $\chi_{rs}^{2}$ are in the same Rouquier block of $(\mathcal{H}_d)_\phi$ for all $0 \leq r<s<d$.

The hyperplane $C_k=C_l$ is $\mathfrak{p}$-essential for the linear characters
$$ \chi_{ijk},\chi_{ijl}\, \textrm{ for all } 0 \leq i,j \leq 1,$$
and the characters of degree $2$
$$\chi_{km}^{1},\chi_{km}^{2},\chi_{lm}^{1},\chi_{lm}^{2}\, \textrm{ for all } 0 \leq m <d \textrm{ with } m \notin \{k,l\}.$$
Due to Proposition $\ref{aA}$, the  formulas of Proposition $\ref{values of aA}$ imply that 
\begin{itemize}
\item the character $\chi_{ijk}$ $(0 \leq i,j\leq 1)$ can be in the same block only with the character $\chi_{ijl}$,
\item the character $\chi_{km}^{1}$ ($ 0 \leq m <d \textrm{ and } m \notin \{k,l\}$) can be in the same block only with the characters $\chi_{km}^{2},\chi_{lm}^{1},\chi_{lm}^{2}$.
\end{itemize}
Let $m \in \{0,1,\ldots,d-1\}\setminus \{k,l\}$. We have that the characters $\chi_{km}^{1}$ and $\chi_{km}^{2}$ are in the same Rouquier block of $(\mathcal{H}_d)_\phi$. The same holds for the characters $\chi_{lm}^{1}$ and $\chi_{lm}^{2}$. Therefore, in order to obtain the desired result, it is enough to show that 
$\{\chi_{km}^{1},\chi_{km}^{2}\}$ isn't a Rouquier block of $(\mathcal{H}_d)_\phi$.

Following \cite{Ma2}, Table $3.10$, there exists an element $T_1$ of $\mathcal{H}_d$ such that
$$\chi_{km}^1(T_1)=\chi_{km}^2(T_1)=x_0+x_1.$$
Let $\mathcal{O}$ be the Rouquier ring of $K$. Suppose that $\{\chi_{km}^{1},\chi_{km}^{2}\}$ is a block of $\mathcal{O}_{\mathfrak{p}\mathcal{O}}(\mathcal{H}_d)_\phi$. Then, by Corollary $\ref{what we are searching}$, we must have
$$\frac{\phi(\chi_{km}^1(T_1))}{\phi(s_{\chi_{km}^1})}+\frac{\phi(\chi_{km}^2(T_1))}{\phi(s_{\chi_{km}^2})}= \phi(x_0+x_1) \cdot \left(\frac{1}{\phi(s_{\chi_{km}^{1}})}+\frac{1}{\phi(s_{\chi_{km}^{2}})}\right) \in \mathcal{O}_{\mathfrak{p}\mathcal{O}}.$$
Since $\phi$ is associated with the hyperplane $C_k=C_l$, we have that $$\phi(x_0+x_1) \notin \mathfrak{p}\mathcal{O}$$ and thus we obtain that
$$\frac{1}{\phi(s_{\chi_{km}^{1}})}+\frac{1}{\phi(s_{\chi_{km}^{2}})} \in \mathcal{O}_{\mathfrak{p}\mathcal{O}}.$$
Using the formulas of Theorem $\ref{nice Schur}$, we can easily calculate that the above element doesn't belong to $\mathcal{O}_{\mathfrak{p}\mathcal{O}}.$

\item Let $\phi$ be a cyclotomic specialization associated with the $\mathfrak{p}$-essential hyperplane $A_i-A_{1-i}+B_j-B_{1-j}+C_k-C_l=0$, where $0 \leq i,j\leq 1$ and  $0 \leq k<l<d$. This hyperplane is $\mathfrak{p}$-essential for the following characters:
$$\chi_{ijk},\chi_{1-i,1-j,l} \textrm{ and  } \chi_{kl}^{1}\textrm{ or }\chi_{kl}^{2}.$$
Let $\mathcal{O}$ be the Rouquier ring of $K$. If the hyperplane is essential for only  three characters, then, due to Proposition $\ref{Malle-Rouquier}$, these three characters are in the same block of $\mathcal{O}_{\mathfrak{p}\mathcal{O}}(\mathcal{H}_d)_\phi$. Otherwise, using the same argument as in the previous case, we can prove that all four characters are in the same block of $\mathcal{O}_{\mathfrak{p}\mathcal{O}}(\mathcal{H}_d)_\phi$.
Now, by Proposition $\ref{Rouquier blocks and central characters}$, the Rouquier blocks of $(\mathcal{H}_d)_\phi$ are unions of the blocks of
$\mathcal{O}_{\mathfrak{p}\mathcal{O}}(\mathcal{H}_d)_\phi$ and $\mathcal{O}_{\mathfrak{I}\mathcal{O}}(\mathcal{H}_d)_\phi$. Therefore, the non-trivial Rouquier blocks of $(\mathcal{H}_d)_\phi$ are
$$\{\chi_{ijk},\chi_{1-i,1-j,l} ,\chi_{kl}^{1},\chi_{kl}^{2}\},$$
$$\{\chi_{rs}^{1},\chi_{rs}^{2}\}\, \textrm{ for all } 0 \leq r<s<d \textrm{ with } (r,s)  \neq (k,l).$$}
\end{enumerate}
\end{apod}

We are now going to prove the following desired result about the functions $a$ and $A$:

\begin{proposition}\label{final a}
Let $\phi: x_i \mapsto (-1)^iq^{a_i}, \,\,
 y_j \mapsto (-1)^jq^{b_j}, \,\,
 z_k \mapsto \zeta_d^kq^{c_k}$
be a cyclotomic specialization for $\mathcal{H}_d$. If the irreducible characters $\chi_\phi$ and $\psi_\phi$ belong to the same Rouquier block of $(\mathcal{H}_{d})_\phi$, then 
$$a_{\chi_\phi}=a_{\psi_\phi} \,\textrm{ and }\, A_{\chi_\phi}=A_{\psi_\phi} .$$
\end{proposition}
\begin{apod}{Thanks to Proposition $\ref{explain AllBlocks}$, it suffices to show that 
the valuations $a_{\chi_\phi}$ and the degrees $A_{\chi_\phi}$ of the Schur elements are constant on the Rouquier blocks associated with an essential hyperplane $H$ (resp. no essential hyperplane), when the integers $a_i,b_j,c_k$ belong to the hyperplane $H$ (resp. no essential hyperplane).

First, due to the description of the Schur elements by Theorem $\ref{nice Schur}$ and the formulas of Proposition $\ref{Aa formula}$, we can deduce that the Schur elements of the characters $\chi_{kl}^1$ and $\chi_{kl}^2$ $(0 \leq k<l<d)$ have the same valuation and the same degree for any cyclotomic specialization $\phi$.

For the same reasons, we have that
\begin{itemize}
\item if $a_0=a_1$, then 
\begin{center}
$a_{\chi_{0jk}}=a_{\chi_{1jk}}$ and $A_{\chi_{0jk}}=A_{\chi_{1jk}}$ for all 
$0\leq j \leq 1,\, 0 \leq k<d$,
\end{center} 
\item if $b_0=b_1$, then 
\begin{center}
$a_{\chi_{i0k}}=a_{\chi_{i1k}}$ and $A_{\chi_{i0k}}=A_{\chi_{i1k}}$ for all 
$0\leq i \leq 1,\, 0 \leq k<d$,
\end{center}
\item if $c_k=c_l$\, $(0 \leq k<l<d)$, then
\begin{center}
$a_{\chi_{ijk}}=a_{\chi_{ijl}}$ and $A_{\chi_{ijk}}=A_{\chi_{ijl}}$ for all 
$0\leq i,j \leq 1$,\\$ $\\
$a_{\chi_{km}^{1,2}}=a_{\chi_{lm}^{1,2}}$ and $A_{\chi_{km}^{1,2}}=A_{\chi_{lm}^{1,2}}$ for all 
$ m \in \{0,1,\ldots,d-1\}\setminus \{k,l\}$.
\end{center}
\end{itemize}

Now let us suppose that $a_i-a_{1-i}+b_j-b_{1-j}+c_k-c_l=0$, with $i,j \in \{0,1\}$,  $k,l \in \{0,1,\ldots,d-1\}$ and $k< l$. We have to show that 
\begin{center}
$a_{\chi_{ijk}}=a_{\chi_{1-i,1-j,l}}=a_{\chi_{kl}^{1,2}}$ and 
$A_{\chi_{ijk}}=A_{\chi_{1-i,1-j,l}}=A_{\chi_{kl}^{1,2}}$.
\end{center}
Due to Proposition \ref{aA}, it suffices to show that
\begin{center}
$a_{\chi_{ijk}}=a_{\chi_{1-i,1-j,l}}=a_{\chi_{kl}^{1,2}}$.
\end{center}
Using the notations of Proposition $\ref{Aa formula}$, Theorem $\ref{nice Schur}$ implies that
\begin{center}
$a_{\chi_{ijk}}=(a_i-a_{1-i})^- + (b_j-b_{1-j})^- +$\\ $\sum_{m=0,\,m\neq k}^{d-1} [(c_k-c_m)^- +
(a_i-a_{1-i}+b_j-b_{1-j}+c_k-c_m)^-],$
\end{center}
\begin{center}
$a_{\chi_{1-i,1-j,l}}=(a_{1-i}-a_{i})^- + (b_{1-j}-b_{j})^- +$\\ 
$\sum_{m=0,\,m\neq l}^{d-1} [(c_l-c_m)^- +
(a_{1-i}-a_i+b_{1-j}-b_{j}+c_l-c_m)^-],$
\end{center}
\begin{center}
$a_{\chi_{kl}^{1,2}}=\sum_{m=0,\,m\neq k,l}^{d-1}[( c_k-c_m)^-+(c_l-c_m)^-]+$ 
$(1/2)\cdot\sum_{h=0}^1 [(a_h-a_{1-h}+b_h-b_{1-h}+c_k-c_l)^- +(a_h-a_{1-h}+b_{1-h}-b_{h}+c_l-c_k)^-].$
\end{center}
Since $a_i-a_{1-i}+b_j-b_{1-j}+c_k-c_l=0$, the above relations give
\begin{center}
$a_{\chi_{ijk}}=(a_i-a_{1-i})^- + (b_j-b_{1-j})^- + \sum_{m=0,\,m\neq k}^{d-1} [(c_k-c_m)^- +
(c_l-c_m)^-],$
\end{center}
\begin{center}
$a_{\chi_{1-i,1-j,l}}=(a_{1-i}-a_{i})^- + (b_{1-j}-b_{j})^- + 
\sum_{m=0,\,m\neq l}^{d-1} [(c_l-c_m)^- +(c_k-c_m)^-],$
\end{center}
\begin{center}
$a_{\chi_{kl}^{1,2}}=
\sum_{m=0,\,m\neq k,l}^{d-1}[( c_k-c_m)^-+(c_l-c_m)^-]+D$,
\end{center}
 where
$$D:=
\left\{
\begin{array}{ll} 
(a_i-a_{1-i})^-+(b_j-b_{1-j})^-+(c_k-c_l)^-,\textrm{ if }i=j,\\
(a_{1-i}-a_i)^-+(b_{1-j}-b_j)^-+(c_l-c_k)^-, \textrm{ if }i \neq j.
\end{array} \right.$$
Obviously, if $i=j$, then
$a_{\chi_{kl}^{1,2}}=a_{\chi_{ijk}}$ and if $i \neq j$, then
$a_{\chi_{kl}^{1,2}}=a_{\chi_{1-i,1-j,l}}$.
Therefore, it is enough to show that 
$a_{\chi_{ijk}}=a_{\chi_{1-i,1-j,l}}$, \ie that
$$(a_i-a_{1-i})^-+(b_j-b_{1-j})^-+(c_k-c_l)^-=(a_{1-i}-a_i)^-+(b_{1-j}-b_j)^-+(c_l-c_k)^-.$$
Since $n^--(-n)^-=n$, for all $n \in \mathbb{Z}$ and
 $a_i-a_{1-i}+b_j-b_{1-j}+c_k-c_l=0$, the above equality holds.}
\end{apod}

\subsection{Rouquier blocks for $G(2pd,2p,2)$}

 Let $p,d \geq 1$. We denote by $\mathcal{H}_{2pd,2p,2}$ the generic
 Hecke algebra of $G(2pd,2p,2)$  generated over the Laurent polynomial ring in $d+4$ indeterminates  
$$\mathbb{Z}[X_0,X_0^{-1},X_1,X_1^{-1},Y_0,Y_0^{-1},Y_1,Y_1^{-1},Z_0,Z_0^{-1},Z_1,Z_1^{-1}\ldots,Z_{d-1},Z_{d-1}^{-1}],$$
by the elements $S,T,U$ satisfying the relations
\begin{itemize}
\item $(S-X_0)(S-X_1)=(T-Y_0)(T-Y_1)=(U-Z_0)(U-Z_1)\cdots(U-Z_{d-1})=0$.
\item $STU=UST$, \,$TUS(TS)^{p-1}=U(ST)^{p}.$ 
\end{itemize}

Let $$\vartheta : \left\{ 
\begin{array}{ll} 
X_i \mapsto (-1)^i q^{a_i} &(0 \leq i \leq 1),\\
Y_j \mapsto (-1)^j q^{b_j} &(0 \leq j \leq 1),\\
Z_k \mapsto \zeta_{d}^{k} q^{c_k} &(0 \leq k <d).  
\end{array} \right. 
$$
be a cyclotomic specialization for  $\mathcal{H}_{2pd,2p,2}$. In order to determine the Rouquier blocks of 
$(\mathcal{H}_{2pd,2p,2})_\vartheta$, we might as well consider the cyclotomic specialization
$$\phi : \left\{ 
\begin{array}{ll} 
X_i \mapsto (-1)^i q^{pa_i} &(0 \leq i \leq 1),\\
Y_j \mapsto (-1)^j q^{pb_j} &(0 \leq j \leq 1),\\
Z_k \mapsto \zeta_{d}^{k} q^{pc_k} &(0 \leq k <d).  
\end{array} \right. 
$$
Since the integers $\{a_i,b_j,c_k\}$ and $\{pa_i,pb_j,pc_k\}$ belong to the same essential hyperplanes for $G(2pd,2p,2)$, Proposition $\ref{explain AllBlocks}$ implies that the Rouquier blocks of 
$(\mathcal{H}_{2pd,2p,2})_\vartheta$ coincide with the Rouquier blocks of 
$(\mathcal{H}_{2pd,2p,2})_\phi$.

We now consider the generic Hecke algebra $\mathcal{H}_{pd}$ of $G(2pd,2,2)$ generated over the ring
$$\mathbb{Z}[x_0,x_0^{-1},x_1,x_1^{-1},y_0,y_0^{-1},y_1,y_1^{-1},z_0,z_0^{-1},z_1,z_1^{-1}\ldots,z_{pd-1},z_{pd-1}^{-1}]$$
by the elements $s,t,u$ satisfying the relations described in the definition of section $4.2$. Let 
$$\phi' : \left\{ 
\begin{array}{ll} 
x_i \mapsto (-1)^i q^{pa_i} &(0 \leq i \leq 1),\\
y_j \mapsto (-1)^j q^{pb_j} &(0 \leq j \leq 1),\\
z_k \mapsto \zeta_{pd}^{k} q^{e_k} &(0 \leq k <pd, e_k: = c_{k \,\mathrm{mod}\,d}).  
\end{array} \right. 
$$
be the ``corresponding'' cyclotomic specialization for  $\mathcal{H}_{pd}$.
Set $\mathcal{H}:=(\mathcal{H}_{pd})_{\phi'}$ and let $\bar{\mathcal{H}}$ be the subalgebra of $\mathcal{H}$ generated by
$s,t$ and $u^p$.
We have
$$(s-q^{pa_0})(s+q^{pa_1})= (t-q^{pb_0})(t+q^{pb_1})=\prod_{k=0}^{d-1}(u^p-\zeta_d^kq^{pc_k})=0.$$
Then (as stated in \cite{Ma2}, Proposition $3.9$)
\cite{Ariki}, Proposition 1.16 implies that the algebra $(\mathcal{H}_{2pd,2p,2})_\phi$ is isomorphic to the algebra $\bar{\mathcal{H}}$ via the morphism
$$S \mapsto  s,\, T \mapsto t,\, U \mapsto u^p.$$

Under the assumptions $\ref{ypo}$, 
the algebra $\mathcal{H}$ is of rank $(2pd)^2$, whereas the algebra 
 $\bar{\mathcal{H}}$ is of rank $(2pd)^2/p$. It is immediate that
 
 \begin{proposition}\label{we can apply clifford}
 The algebra  $\mathcal{H}$ is a free  $\bar{\mathcal{H}}$ -module with basis $\{1,u,\ldots,u^{p-1}\}$, i.e.,
 $$\mathcal{H}= \bar{\mathcal{H}}  \oplus  u\bar{\mathcal{H}} \oplus \cdots \oplus u^{p-1}\bar{\mathcal{H}}.$$  
\end{proposition}

Again under the assumptions $\ref{ypo}$, the algebra $\mathcal{H}$ is symmetric and $\bar{\mathcal{H}}$ is a symmetric subalgebra of $\mathcal{H}$. In particular, following Definition $\ref{symmetric algebra of a finite group}$, $\mathcal{H}$ is the twisted symmetric algebra of the cyclic group of order $p$ over $\bar{\mathcal{H}}$ (since $u$ is a unit in $\mathcal{H}$). Therefore, we can apply Proposition $\ref{1.45}$ and obtain (using the notations of section $1.3$) the following.

\begin{proposition}\label{1}
If $G$ is the cyclic group of order $p$ and $K:=\mathbb{Q}(\zeta_{2pd})$, then the block-idempotents of $(Z\mathcal{R}_K(q)\bar{\mathcal{H}})^G$ coincide with the block-idempotents of $(Z\mathcal{R}_K(q)\mathcal{H})^{G^\vee}$, where $\mathcal{R}_K(q)$ is the Rouquier ring of $K$.
\end{proposition}

The action of the cyclic group $G^\vee$ of order $p$ on $\mathrm{Irr}(K(q)\mathcal{H})$ corresponds to the action 
$$\chi_{i,j,k} \mapsto \chi_{i,j,k+d }\,\, (0\leq i,j\leq 1) \,(0 \leq k <pd),$$
$$\chi_{k,l}^{1,2} \mapsto \chi_{k+d ,l+d }^{1,2}\,\,
(0 \leq k <l<pd),$$
where all the indexes are considered  $\mathrm{mod}\,pd$.
%
With the help of the following lemma, we will show that the Rouquier blocks of 
$\mathcal{H}$ are stable under the action of $G^\vee$. Here the results of Theorem $\ref{yes proof}$ are going to be used as definitions.

\begin{lemma}\label{three hyperplanes}
Let $k_1$, $k_2$ and $k_3$ be three distinct elements
of  $\{0,1,\ldots,pd-1\}$. If  the blocks of $\mathcal{R}_K(q)\mathcal{H}$ are unions of the Rouquier blocks associated with the (not necessarily essential) hyperplanes $C_{k_1}=C_{k_2}$ and $C_{k_2}=C_{k_3}$, then they are also unions of the Rouquier blocks associated with the (not necessarily essential) hyperplane
 $C_{k_1}=C_{k_3}$.
\end{lemma}
\begin{apod}{We only need to show that
\begin{enumerate}[(a)]
\item the characters $\chi_{i,j,k_1}$ and $\chi_{i,j,k_3}$ are in the same block of $\mathcal{R}_K(q)\mathcal{H}$ for all $0 \leq i,j \leq 1,$
\item the characters $\chi_{k_1,m}^{1,2}$ and $\chi_{k_3,m}^{1,2}$
are in the same block of $\mathcal{R}_K(q)\mathcal{H}$ for all
 $0 \leq m <pd$ with  $m \notin \{k_1,k_3\}.$
\end{enumerate}
Since the blocks of $\mathcal{R}_K(q)\mathcal{H}$ are unions of the Rouquier blocks associated with the  hyperplanes $C_{k_1}=C_{k_2}$ and $C_{k_2}=C_{k_3}$, Theorem $\ref{yes proof}$ yields that
\begin{enumerate}[(1)]
\item the characters $\chi_{i,j,k_1}$ and $\chi_{i,j,k_2}$ are in the same block of $\mathcal{R}_K(q)\mathcal{H}$ for all $0 \leq i,j \leq 1,$
\item the characters $\chi_{i,j,k_2}$ and $\chi_{i,j,k_3}$ are in the same block of $\mathcal{R}_K(q)\mathcal{H}$ for all $0 \leq i,j \leq 1,$
\item the characters $\chi_{k_1,m}^{1,2}$ and $\chi_{k_2,m}^{1,2}$
are in the same block of $\mathcal{R}_K(q)\mathcal{H}$ for all
 $0 \leq m <pd$ with  $m \notin \{k_1,k_2\},$
\item the characters $\chi_{k_2,m}^{1,2}$ and $\chi_{k_3,m}^{1,2}$
are in the same block of $\mathcal{R}_K(q)\mathcal{H}$ for all
 $0 \leq m <pd$ with  $m \notin \{k_2,k_3\}.$
\end{enumerate}
We immediately deduce (a) for all $0 \leq i,j \leq 1$ and (b) for all $0 \leq m <pd$ with $m \notin \{k_1,k_2,k_3\}.$
Finally, (3) implies that the characters $\chi_{k_1,k_3}^{1,2}$ and $\chi_{k_2,k_3}^{1,2}$
are in the same block of $\mathcal{R}_K(q)\mathcal{H}$, whereas by (4), $\chi_{k_1,k_2}^{1,2}$ and $\chi_{k_1,k_3}^{1,2}$ are also in the same block of $\mathcal{R}_K(q)\mathcal{H}$.
Thus, the characters $\chi_{k_1,k_2}^{1,2}$ and $\chi_{k_2,k_3}^{1,2}$ belong to the same Rouquier block of $\mathcal{H}$.}
\end{apod}

\begin{theorem}\label{2}
The blocks of 
$\mathcal{R}_K(q)\mathcal{H}$ are stable under the action of $G^\vee$.
\end{theorem}
\begin{apod}{Following Proposition $\ref{explain AllBlocks}$, the Rouquier blocks of $\mathcal{H}$ are unions of the Rouquier blocks associated with all the essential hyperplanes of the form
$$C_{h+md}=C_{h+nd} \,\,(0 \leq h < d,\,0 \leq m < n < p).$$
Recall that the hyperplane
$C_{h+md}=C_{h+nd} $ is actually essential for $G(2pd,2,2)$ if and only if
the element $\zeta_{pd}^{h+md}-\zeta_{pd}^{h+nd}$ belongs to a prime ideal of
$\mathbb{Z}[\zeta_{2pd}]$, \ie if and only if the element $\zeta_{p}^{m}-\zeta_{p}^{n}$ belongs to a prime ideal of $\mathbb{Z}[\zeta_{2pd}]$.

Suppose that $p=p_1^{t_1}p_2^{t_2}\cdots p_r^{t_r}$, where the $p_i$ are distinct prime numbers. For $s \in \{1,2,\ldots,r\}$, we set
$h_s:= p/p_s^{t_s}.$ Then $\mathrm{gcd}(h_s)=1$ and by Bezout's theorem, there exist integers $(g_s)_{1 \leq s \leq r}$ such that $\sum_{s=1}^rg_sh_s=1$.
The element $1-\zeta_{p}^{g_sh_s}$ belongs to all the prime ideals of $\mathbb{Z}[\zeta_{2pd}]$ lying over the prime number $p_s$.  Let $h \in \{0,1,\ldots,d-1\}$ and $m \in \{0,1,\ldots,p-2\}$ and set
 $$l_0:=m \textrm{ and } l_s:=(l_{s-1}+g_sh_s) \,\,\mathrm{ mod }\,p, \textrm{ for all } s\,(1 \leq s \leq r).$$
We have that the element $\zeta_{p}^{l_{s-1}}-\zeta_{p}^{l_{s}}=\zeta_{p}^{l_{s-1}}(1-\zeta_{p}^{g_{s}h_s})$ belongs to all the prime ideals of $\mathbb{Z}[\zeta_{2pd}]$ lying over the prime number $p_s$.
Therefore, the hyperplane $C_{h+l_{s-1}d}=C_{h+l_sd}$ is essential for $G(2pd,2,2)$ for all $s$ $(1 \leq s \leq r)$.
Since $l_0=m$ and $l_r=m+1$, Lemma $\ref{three hyperplanes}$ implies that the Rouquier blocks of $\mathcal{H}$ are unions of the Rouquier blocks associated with the (not necessarily essential) hyperplane
$$C_{h+md}=C_{h+(m+1)d},$$
following their description by Theorem $\ref{yes proof}$. Since this holds for all
$m\,(0 \leq m \leq p-2)$, Lemma $\ref{three hyperplanes}$ again implies that 
the Rouquier blocks of $\mathcal{H}$ are unions of the Rouquier blocks associated with  all the hyperplanes of the form 
$$C_{h+md}=C_{h+nd} \,\,(0 \leq m < n < p),$$
for all $h\,(0 \leq h < d)$. We deduce that
\begin{enumerate}[(1)]
\item the characters $(\chi_{i,j,h+md})_{ 0 \leq m < p}\,$ belong to the same block of  $\mathcal{R}_K(q)\mathcal{H}$, for all $0\leq i,j\leq 1$ and $0\leq h <d,$
\item the characters $(\chi_{h+md,h+nd}^{1,2})_{ 0 \leq m<n< p}\,$ belong to the same block of  $\mathcal{R}_K(q)\mathcal{H}$, for all  $0\leq h <d,$
\item the characters $(\chi_{h+md,h'+nd}^{1,2})_{ 0 \leq m,n< p}\,$ belong to the same block of  $\mathcal{R}_K(q)\mathcal{H}$, for all  $0\leq h<h' <d$.
\end{enumerate}
Hence, the blocks of $\mathcal{R}_K(q)\mathcal{H}$ are stable under the action of $G^\vee$.
} 
\end{apod}

Following Theorem $\ref{2}$, Proposition $\ref{1}$ now gives:

\begin{corollary}\label{3}
If $G$ is the cyclic group of order $p$ and $K:=\mathbb{Q}(\zeta_{2pd})$, then the block-idempotents of $(Z\mathcal{R}_K(q)\bar{\mathcal{H}})^G$ coincide with the block-idempotents of $\mathcal{R}_K(q)\mathcal{H}$.
\end{corollary}

Now, let $\bar{\chi} \in \mathrm{Irr}(K(q)\bar{\mathcal{H}})$. Using the notations of Proposition $\ref{1.42}$, we have that
$|\Omega||\bar{\Omega}|=p.$
Since $|\Omega|=p$, we obtain that $|\bar{\Omega}|=1$
and thus $e(\bar{\chi})$ is fixed by the action of $G$. Therefore,
the block-idempotents of $\mathcal{R}_K(q)\bar{\mathcal{H}}$ are also fixed by the action of $G$. Consequently, we obtain the following.

\begin{proposition}\label{4}
The block-idempotents of $\mathcal{R}_K(q)\bar{\mathcal{H}}$ coincide with the block-idempotents of $\mathcal{R}_K(q)\mathcal{H}$.
\end{proposition}

 Thanks to the above result, in order to determine the Rouquier blocks of $\bar{\mathcal{H}}$, it suffices to calculate the Rouquier blocks of $\mathcal{H}$ and restrict all the characters to $\bar{\mathcal{H}}$. The Rouquier blocks of $\mathcal{H}$ can be obtained with the use of Theorem $\ref{yes proof}$.\\

Now,
\begin{itemize}
\item the description of the Rouquier blocks of $\bar{\mathcal{H}}$ by Proposition $\ref{3}$,
\item the relation between the Schur elements of $\bar{\mathcal{H}}$
and the Schur elements of $\mathcal{H}$ given by Proposition $\ref{1.42}$
\item and the invariance of the integers $a_\chi$ and $A_\chi$ on the Rouquier blocks of $\mathcal{H}$, resulting from Proposition $\ref{final a}$, imply that
\end{itemize}

\begin{proposition}\label{fin}
The valuations $a_{\bar{\chi}}$ and the degrees $A_{\bar{\chi}}$ of the Schur elements are constant on the Rouquier blocks of $\bar{\mathcal{H}}$.
\end{proposition}

\end{document}